\documentclass[11pt]{article}
\usepackage{amsmath, amsthm, amssymb, }
\usepackage{float,epsfig}
\usepackage{subcaption, multirow}
\usepackage{tkz-graph}
\usepackage{amsthm}
\usepackage{color,graphicx,float,lscape,enumerate}
\usepackage[colorlinks = true]{hyperref}
\usepackage{cleveref}
\usepackage{cite}

\newtheorem{theorem}{\bf Theorem}[section]

\newtheorem{definition}[theorem]{\bf Definition}

\newtheorem{remark}[theorem]{\bf Remark}
\newtheorem{lemma}[theorem]{\bf Lemma}

\newcommand{\ba}{\begin{array}}
\newcommand{\ea}{\end{array}}

\newcommand{\vone}{\vskip 2ex}

\newcommand{\be}{\begin{equation}}
\newcommand{\ee}{\end{equation}}
\newcommand{\beano}{\begin{eqnarray*}}
\newcommand{\eeano}{\end{eqnarray*}}

\def\bmatrix#1{\left[ \begin{matrix} #1 \end{matrix} \right]}

\def \R{{\mathbb R}}

\def \M{{\mathrm M}}

\def \pf{{\bf Proof: }}

\def \T{{\mathrm T}}

\def \lam{{\lambda}}

\title{Corona product of signed graphs and its application to signed network modelling}
\author{ Bibhas Adhikari  \,\, Amrik Singh \,\, Sandeep Kumar Yadav}
\date{}
\begin{document}
\maketitle

{\small \noindent{\bf Abstract.} The notion of corona of two graphs was introduced by Frucht and Harary in 1970. In this paper we generalize their definition of corona product of two graphs and introduce corona product of two signed graphs by utilizing the framework of marked graphs, which was introduced by Beineke and Harary in 1978. We study structural and spectral properties of corona product of signed graphs. Further we define signed corona graphs by considering corona product of a fixed small signed graph  with itself iteratively, and we call the small graph as the seed graph for the corresponding corona product graphs. Signed corona graphs can be employed as a signed network generative model for large growing signed networks. We study structural properties of corona graphs that include statistics of signed links, all types of signed triads and degree distribution. Besides we analyze algebraic conflict of signed corona graphs generated by specially structured seed graphs. Finally we show that a suitable choice of a seed graph can produce corona graphs which preserve properties of real signed networks.

\vone \noindent{\bf Keywords.} signed graphs, structural balance, corona product, algebraic conflict


\section{Introduction}\label{sec:1}

Signed networks represent a framework to deal with binary relationship between nodes in a network that has two contradictory possibilities. For example, like and dislike, love and hate, trust and distrust are considered as measures of relationships between people, whereas alliance and antagonism between two countries can be considered as contradictory binary international relationships.  Signed network model to represent social systems was first introduced by Harary and Cartwright in 1956 to generalize the theory of \textit{balanced state} of a social system developed by Heider in 1946 \cite{cartwright1956structural}. Heider rationalized his theory of balanced state by considering possible relationships in a system of three entities \cite{heider1946attitudes}. On the other hand, real world data can be put into the signed network setup to determine salient features of the data. For instance, in the trust network of Epinions, users establish binary relations to each other that reflect trust and distrust \cite{guha2004propagation,leskovec2010signed}. Mathematically, a signed graph is a graph in which some edges are designated with positive sign reflecting positive relationship between the constituent pair of nodes and other edges are assigned negative sign which represents negative relationship. Thus signed graph is an ordered tuple $G=(V, E, \sigma)$ where $V$ denotes the set of nodes, $E\subseteq V\times V$ the edge set, and $\sigma: E\rightarrow \{+,-\}$ is called the signature function \cite{zaslavsky1982signed}. In this paper we use the terms signed graph and signed network interchangeably.

Structural balance of signed networks has always been the central topic in the study of real world signed networks. A signed network is called balanced or in balanced state if all triads in it are balanced \cite{cartwright1956structural}. A triad is balanced if all three edges are positive, or if two are positive and one negative. Alternatively, it is also established that a signed network is balanced if and only if all its cycles are balanced \cite{harary1953notion}. A signed cycle is called balanced if the number of negative edges in it is even. Most often it is hypothesized that real networks evolve towards balanced state \cite{doreian2013brief}. See also \cite{doreian2004evolution}. However, recently an excellent article including a historical background about the development of the theory of balance demonstrates that this is not the case always \cite{estrada2019rethinking}. If a signed network is not balanced it is called unbalanced. Evidently, it is unfair to call a signed network unbalanced if only a few cycles in it are unbalanced compared to an unbalanced network in which most of the cycles are unbalanced. Thus the concept of degree of balance of a signed network has emerged in literature, and several measures are proposed to estimate it \cite{aref2017measuring,aref2018balance,estrada2014walk,singh2017measuring,aref2019signed,kunegis2010spectral}.    

Several signed network generative models are proposed in literature with a goal to preserve structural or statistical properties of real world signed networks. For example, in \cite{hsieh2012low} a low-rank model based on matrix completion technique is proposed to generate signed  networks which can inherit structural balance property of real networks. In 
\cite{derr2018signed}, the authors propose a parametric model to preserve degree distribution, sign distribution, and balance/unbalanced triad distribution  of signed networks. This model is inspired by the popular Chung-Lu model for generation of unsigned networks. Recently, a model is proposed based on local preferential attachment to generate signed networks that have community structure and high positive clustering coefficient \cite{li2018modeling}. To the best of the knowledge of the authors, the literature lack a deterministic growing network model that can preserve  sign distribution, and balance/unbalanced triad distribution of signed networks. Besides, spectral properties of the existing network generative models are completely unexplored.Thus how the spectral properties relate to the structural properties of the networks generated by these models is unclear. Spectral property such as algebraic conflict \cite{aref2019signed,kunegis2010spectral}, the smallest signed Laplacian eigenvalue of the network can estimate the degree of balance. 


In this paper we first introduce the notion of corona product of signed graphs by generalizing the definition of corona product of unsigned graphs. Let $G$ be an unsigned graph on $n$ nodes, and $H$ an unsigned graph. Then the corona product graph $G\circ H$ of $G$ and $H$ is defined as the graph obtained by taking one copy of $G$ and $n$ copies of $H$, and then joining the $i$th vertex of $G$ to every vertex in the $i$th copy of $H$ \cite{frucht1970corona}. Note that except the existing edges of $G, H$ in $G\circ H$ there will be $nk$ new edges created where $k$ denotes the number of nodes in $H.$ For corona product of signed graphs, the job is to define the signs of these new edges so that when the graphs $G, H$ are unsigned those edge signs will be positive. We define the signs of the new edges utilizing the formalism of marked graph which is a framework introduced by Beineke and Harary \cite{beineke1978consistent,harary1953notion}. We call this a marking scheme on the node set for the definition of corona product. We study structural balance of $G\circ H,$ and provide a necessary and sufficient condition based on the structural properties of $G, H$ such that $G\circ H$ becomes balanced. Let $T_i$ denote a triad in a signed network with $i$ number of negative edges in it, $i=0,1,2,3.$ Obviously, the distribution of the numbers of $T_i$s influence degree of balance of a signed network. Some of the earliest measures of degree of balance are defined in terms of the number of $T_i$s. For example, the ratio of the number of signed to unsigned triads in a signed network is considered as the degree of balance.

\begin{figure}[H]
			\centering	
			\begin{subfigure}[b]{0.2\textwidth}
				\centering
				\begin{tikzpicture}
				\draw [fill] (0, 0) circle [radius=0.1];
				\draw [fill] (1, 0) circle [radius=0.1];
				\draw (0,0) --(1,0);
				\draw [fill] (0.5, 1) circle [radius=0.1];
				\draw (0,0) -- (0.5,1);
				\draw (1,0) -- (0.5,1);
				\end{tikzpicture}
				\caption{$T_0$}
			\end{subfigure}
			\begin{subfigure}[b]{0.2\textwidth}
				\centering
				\begin{tikzpicture}
				\draw [fill] (0, 0) circle [radius=0.1];
				\draw [fill] (1, 0) circle [radius=0.1];
				\draw [dashed] (0,0) --(1,0);
				\draw [fill] (0.5, 1) circle [radius=0.1];
				\draw (0,0) -- (0.5,1);
				\draw (1,0) -- (0.5,1);
				\end{tikzpicture}
				\caption{$T_1$}
			\end{subfigure}		
			\begin{subfigure}[b]{0.2\textwidth}
				\centering
				\begin{tikzpicture}
				\draw [fill] (0, 0) circle [radius=0.1];
				\draw [fill] (1, 0) circle [radius=0.1];
				\draw (0,0) --(1,0);
				\draw [fill] (0.5, 1) circle [radius=0.1];
				\draw [dashed] (0,0) -- (0.5,1);
				\draw [dashed] (1,0) -- (0.5,1);
				\end{tikzpicture}
				\caption{$T_2$}
			\end{subfigure}			
			\begin{subfigure}[b]{0.20\textwidth}
				\centering
				\begin{tikzpicture}
				\draw [fill] (0, 0) circle [radius=0.1];
				\draw [fill] (1, 0) circle [radius=0.1];
				\draw [dashed] (0,0) --(1,0);
				\draw [fill] (0.5, 1) circle [radius=0.1];
				\draw [dashed] (0,0) -- (0.5,1);
				\draw [dashed] (1,0) -- (0.5,1);
				\end{tikzpicture}
				\caption{$T_3$}
			\end{subfigure}
			\caption{Possible triads in an (undirected) signed network. Triads (a) and (c) are balanced, and triads (b) and (d) are unbalanced. Solid and dashed edges represent positive and negative edges, respectively.}
			\label{fig:triads}
		\end{figure}
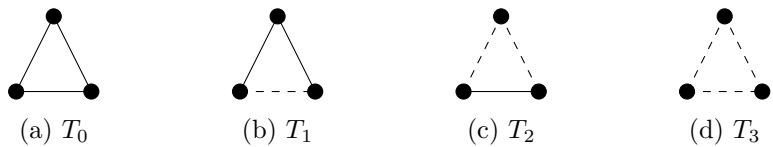	

The positive and negative degree of a node $u$ in a signed graph are defined by number of positive and negative edges incidental to the node $u$ denoted by $d^+(u)$ and $d^-(u)$, respectively. The total degree and net-degree of $u$ are given by $d^\pm(u) =d^+(u) + d^-(u)$ and  $d(u) =d^+(u) - d^-(u)$, respectively. A signed graph is called net-regular if every node of it has the same net-degree $d$ \cite{netregular}. Note that net-regular signed graphs are signed counterpart of unsigned regular graphs.  The adjcency matrix $A=[a_{ij}]$ corresponding to a signed graph $G=(V,E,\sigma)$ is defined by $a_{ij}=1$ if $(i,j)\in E$ and $\sigma((i,j))=+,$ $a_{ij}=-1$ if $\sigma((i,j))=-,$ and $a_{ij}=0$ otherwise. The signed Laplacian matrix of $G$ on $n$ nodes is defined as $L=D-A$ where $D$ is a diagonal matrix of order $n$ whose $i$th diagonal entry is the total degree of $i$th node of $G$ \cite{kunegis2010spectral}. The signless Laplacian matrix of $G$ is defined by $Q=D+A.$ The eigenvalues of $A, L,$ and $Q$ are called adjacency, signed Laplacian and signless Laplacian eigenvalues of $G.$ The eigenvalues of these matrices inherit structural balance properties of the corresponding signed graph. For example, a signed graph $G$ is balanced if and only if adjacency eigenvalues of $G$ are same as the adjacency eigenvalues of the underlying unsigned graph corresponding to $G$ \cite{acharya1980spectral}. The least signed Laplacian eigenvalue is zero if and only if the corresponding graph is balanced \cite{kunegis2010spectral}. For an unbalanced graph, the least signed Laplacian eigenvalues is considered as a measure of degree of balance of the graph and it is called algebraic conflict of the graph \cite{kunegis2010spectral,aref2019signed}. We derive adjacency eigenvalues of $G\circ H$ in terms of adjacency eigenvalues of the signed graphs $G$ and $H$ respectively, when $H$ is net-regular. We determine signed Laplacian eigenvalues of $G\circ H$ from which the formula of signless Laplacian eigenvalues follows after a marginal modification.

Utilizing the corona product framework for unsigned graphs Sharma et al. \cite{sharma2015spectra,sharma2017structural} defined corona graphs as a model for generating large networks. See also \cite{lv2015corona}. Given a connected graph $G=G^{(0)}$ the corona graphs are defined as \begin{equation}\label{def:cg} G^{(m)}=G^{(m-1)}\circ G,\,\, m\geq 1.\end{equation} The graph $G$ is called the seed graph for the corona graph $G^{(m)}.$ Adapting this idea and incorporating the definition of corona product introduced in this paper, we define signed corona graphs by considering the seed graph as a signed graph as in equation (\ref{def:cg}). Thus construction of signed corona graphs provides a deterministic growing signed network generative model. Obviously, the structural balance and spectral properties of $G^{(m)}$ depend on the choice of the seed graph and the distribution of signs of the new edges which appear due to the definition of corona product, that is, the marking scheme on the node set of the seed graph. The modelling perspective of preserving properties of real networks using $G^{(m)}$ raises the following question. Can we choose a seed graph $G$ so that desired sign distribution, triad distribution, degree of balance in $G^{(m)}$ be obtained? We show that this is indeed the case. That is, given the statistics of signs of edges and triads of a real signed network it is possible to determine a seed graph such that the corresponding large signed corona graphs provide desired distribution of signed edges and triads approximately. The algebraic conflict of $G^{(m)}$ can be obtained by iterative use of the formula of signed Laplacian eigenvalues derived for corona product of signed graphs discussed above. Besides, we determine degree sequence of $G^{(m)}$ from which the positive and negative degree distribution of $G^{(m)}$ can be computed and compared with real networks. 


Throughout we assume that the seed graph is a simple connected signed graph. The main contributions of this paper are as follows
\begin{itemize}
\item[(1)] We introduce corona product of two signed graphs based on the framework of marked graphs and study its structural properties. \item[(2)] We study spectra and Laplacian spectra of corona product of signed graphs, and provide computable expressions of eigenvectors corresponding to such eigenvalues. \item[(3)] We define signed corona graphs and propose it as a signed network generative model. We provide computable formulae of number of signed edges, triads, and positive and negative degrees of every node in a corona graph. We investigate algebraic conflict of corona graphs, which are generated by specially structured seed graphs. Finally we show that corona graphs can inherit properties of real signed networks. 
\end{itemize}

\section{Corona product of signed graphs}

In this section we introduce corona product $G\circ H$ of a pair of signed graphs $G, H$, and we study the structural balance and spectral properties of $G\circ H.$ Note that this operation (product) is non-commutative, that is, $G\circ H$ need not be equal to $H\circ G.$ First we review the notion of marked graphs as follows.

A graph is called a marked graph if every node of the graph is marked by either a positive or negative sign. Thus a marked graph is a tuple $G=(V, E,\mu)$ where $V$ is the node set, $E$ the edge set and $\mu: V\rightarrow \{+,-\}$ is called the marking function.  An obvious way to construct a signed graph from a marked graph is be defining the sign of an edge of the marked graph as the product of signs of its adjacent vertices \cite{harary1981counting}. On the other hand, a marked graph can be defined from a signed graph $G=(V, E,\sigma)$ by defining the marking of a node $v\in V$ as \begin{equation}\label{Def:markingofanode} \mu(v)=\prod_{e\in E_v} \sigma(e)\end{equation} where $E_v$ is the set of signed edges adjacent at $v.$ Such a method of marking is also known as \textit{canonical marking} \cite{beineke1978consistent,beineke1978consistent2}. There can be multiple ways to define a marking function for a signed graph. We consider two marking functions, canonical marking and \textit{plurality marking} in this paper. We define plurality marking of a node $v$ of a signed graph $G=(V,E,\mu)$ as  \begin{equation}\label{Def:mmarkingofanode}\mu(v)=
\left\{
  \begin{array}{ll}
    +, & \mbox{if} \,\, \hbox{$\max\{d^+(v), d^-(v)\} =d^+(v)$} \\
    -, & \mbox{Otherwise}
  \end{array}
\right.
\end{equation} Hence a node is negatively marked in plurality marking scheme only when $d^-(v)>d^{+}(v).$ Thus from now onward we denote a signed graph as a $4$-tuple $G=(V,E,\sigma,\mu)$ where $\sigma$ and $\mu$ are the signature function and marking function defined on the edge set $E$ and node set $V$ respectively.

Thus we introduce the following definition for corona product of two signed graphs.
\begin{definition}
Let $G_1=(V_1,E_1,\sigma_1,\mu_1)$  and $G_2=(V_2,E_2,\sigma_2,\mu_2)$ be signed graphs on $n$ and $k$ nodes respectively. Then corona product $G_1\circ G_2$ of $G_1, G_2$ is a signed graph by taking one copy of $G_1$ and $n$ copies of $G_2,$ and then forming a signed edge from $i$th node of $G_1$ to every node of the $i$th copy of $G_2$ for all $i.$ The sign of the new edge between $i$th node of $G_1,$ say $u$ and $j$th node in the $i$th copy of $G_2,$ say $v$ is given by $\mu_1(u)\mu_2(v)$ where $\mu_i$ is a marking scheme defined by $\sigma_i, i=1,2.$
\end{definition}

For instance, the corona product $G_1\circ G_2$ of signed graphs $G_1$ and $G_2$ is shown in Figure \ref{CP1}. Note that canonical and plurality marking are same for the graph $G_2.$ For $G_1$ the marking of the nodes $1, 3$ are same for canonical and plurality markings, whereas the canonical and plurality markings of node $2$ are $-$ and $+$ respectively. Thus the choice of the marking function produce different corona product graphs. 
		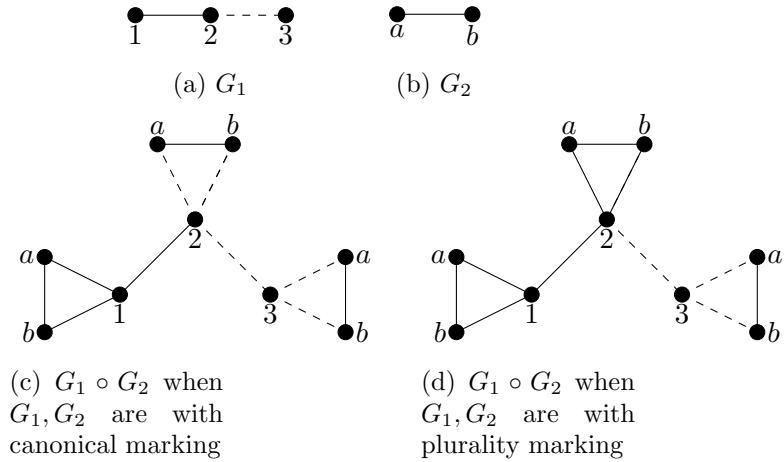
\begin{figure}[H]
			\centering			
			\begin{subfigure}{0.2\textwidth}
				\centering
				\begin{tikzpicture}
				\draw [fill] (0, 0) circle [radius=0.1];
				\node [below] at (0, 0) {$1$};
				\draw [fill] (1, 0) circle [radius=0.1];
				\node [below] at (1, 0) {$2$};
				\draw [fill] (2, 0) circle [radius=0.1];
				\node [below] at (2, 0) {$3$};
				\draw (0,0) --(1,0);
				\draw [dashed] (1, 0) -- (2, 0);
				\end{tikzpicture}
				\caption{$G_1$}
			\end{subfigure}
			\begin{subfigure}{0.2\textwidth}
				\centering
				\begin{tikzpicture}
				\draw [fill] (0, 0) circle [radius=0.1];
				\node [below] at (0, 0) {$a$};
				\draw [fill] (1, 0) circle [radius=0.1];
				\node [below] at (1, 0) {$b$};
				\draw (0,0) --(1,0);
				\end{tikzpicture}
				\caption{$G_2$}
			\end{subfigure} \\
			\begin{subfigure}{0.2\textwidth}
				\centering
				\begin{tikzpicture}
				\draw [fill] (0, -0.5) circle [radius = .1];
				\node [below] at (0, -0.5) {};
				\node [left] at (0, -0.5) {$b$};
				\draw [fill] (1, 0) circle [radius = .1];
				\node [below] at (1, 0) {$1$};
				\node [above] at (1, 0) {};
				\draw (0,-0.5) -- (1, 0);
				\draw [fill] (0, 0.5) circle [radius = .1];
				\node [above] at (0, .5) {};
				\node [left] at (0, .5) {$a$};
				\draw (0,-0.5) -- (0, .5);
				\draw (0,0.5) -- (1, 0);
				\draw [fill] (2, 1) circle [radius = .1];
				\node [below] at (2, 1) {$2$};
				\node [left] at (2, 1) {};
				\draw (2,1) -- (1, 0);
				\draw [fill] (1.5, 2) circle [radius = .1];
				\node [above] at (1.5, 2) {$a$};
				\node [left] at (1.5, 2) {};
				\draw (1.5, 2) -- (2.5, 2);
				\draw [fill] (2.5, 2) circle [radius = .1];
				\node [above] at (2.5, 2) {$b$};
				\node [right] at (2.5, 2) {};
				\draw [dashed] (2, 1) -- (1.5, 2);
				\draw [dashed] (2, 1) -- (2.5, 2);
				\draw [dashed] (2, 1) -- (3, 0);
				\draw [fill] (3, 0) circle [radius = .1];
				\node [below] at (3, 0) {$3$};
				\node [left] at (3, 0) {};
				\draw [dashed] (3, 0) -- (4, -0.5);
				\draw [dashed] (3, 0) -- (4, 0.5);
				\draw [dashed] (2, 1) -- (2.5, 2);
				\draw (4, 0.5) -- (4, -0.5);
				\draw [fill] (4, -0.5) circle [radius = .1];
				\node [below] at (4, -0.5) {};
				\node [right] at (4,-0.5) {$b$};
				\draw [fill] (4, 0.5) circle [radius = .1];
				\node [above] at (4, 0.5) {};
				\node [right] at (4,0.5) {$a$};
				\end{tikzpicture}
				\caption{$G_1 \circ G_2$ when $G_1, G_2$ are with canonical marking}
			\end{subfigure}\hspace{2.5cm}
	\begin{subfigure}{0.2\textwidth}
				\centering
				\begin{tikzpicture}
				\draw [fill] (0, -0.5) circle [radius = .1];
				\node [below] at (0, -0.5) {};
				\node [left] at (0, -0.5) {$b$};
				\draw [fill] (1, 0) circle [radius = .1];
				\node [below] at (1, 0) {$1$};
				\node [above] at (1, 0) {};
				\draw (0,-0.5) -- (1, 0);
				\draw [fill] (0, 0.5) circle [radius = .1];
				\node [above] at (0, .5) {};
				\node [left] at (0, .5) {$a$};
				\draw (0,-0.5) -- (0, .5);
				\draw (0,0.5) -- (1, 0);
				\draw [fill] (2, 1) circle [radius = .1];
				\node [below] at (2, 1) {$2$};
				\node [left] at (2, 1) {};
				\draw (2,1) -- (1, 0);
				\draw [fill] (1.5, 2) circle [radius = .1];
				\node [above] at (1.5, 2) {$a$};
				\node [left] at (1.5, 2) {};
				\draw (1.5, 2) -- (2.5, 2);
				\draw [fill] (2.5, 2) circle [radius = .1];
				\node [above] at (2.5, 2) {$b$};
				\node [right] at (2.5, 2) {};
				\draw (2, 1) -- (1.5, 2);
				\draw (2, 1) -- (2.5, 2);
				\draw [dashed] (2, 1) -- (3, 0);
				\draw [fill] (3, 0) circle [radius = .1];
				\node [below] at (3, 0) {$3$};
				\node [left] at (3, 0) {};
				\draw [dashed] (3, 0) -- (4, -0.5);
				\draw [dashed] (3, 0) -- (4, 0.5);
				\draw [dashed] (2, 1) -- (2.5, 2);
				\draw (4, 0.5) -- (4, -0.5);
				\draw [fill] (4, -0.5) circle [radius = .1];
				\node [below] at (4, -0.5) {};
				\node [right] at (4,-0.5) {$b$};
				\draw [fill] (4, 0.5) circle [radius = .1];
				\node [above] at (4, 0.5) {};
				\node [right] at (4,0.5) {$a$};
				\end{tikzpicture}
				\caption{$G_1 \circ G_2$ when $G_1, G_2$ are with plurality marking}
			\end{subfigure}
			\caption{The corona product of $G_1\circ G_2$ is shown in (c), (d) with canonical and plurality marking functions on $G_i$, $i=1,2$.}
			\label{CP1}
		\end{figure}

Below we discuss structural and spectral properties of the corona product of two graphs.

\subsection{Structural and spectral properties of $G_1\circ G_2$}

We provide the statistics about the number of edges and triads in $G=G_1\circ G_2=(V, E,\sigma,\mu)$ for a given pair of signed graphs  $G_1=(V_1,E_1,\sigma_1,\mu_1)$  and $G_2=(V_2,E_2,\sigma_2,\mu_2)$ as follows. Observe that the number of nodes in $G_1\circ G_2$ is given by $|V_1| + |V_1| \, |V_2|$ where $|V_i|$ denotes the number of nodes in $G_i,$ $i=1,2.$  Let $\mathrm{M}_i^+$ and $\mathrm{M}_i^-$ denote the number of positively and negatively marked nodes in $G_1$ and $G_2$ respectively. The Table \ref{tab:stat} describes the statistics of edges in $G_1\circ G_2.$

\begin{table}[H]
			\centering
			\begin{tabular}{l|c|c|c}
				Edges & $G_1$ & $G_2$ & $G_1\circ G_2$ \\ \hline\hline
				$\#$ of Edges & $|E_1|$ & $|E_2|$ & $|E_1| + |V_1|\,|E_2| + |V_1|\,|V_2|$ \\ 
				$\#$ of $+$ edges & $|E_1^+|$ & $|E_1^-|$ & $|E_1^+| + |V_1| \, |E_2^+| + \M_1^+ \,M_2^+ + M_1^- \, M_2^-$ \\ 
				$\#$ of $-$ edges & $|E_1^-|$ & $|E_2^-|$ & $|E_1^-| + |V_1| \, |E_2^-| + \M_1^+ \,M_2^- + M_1^- \, M_2^+$ \\ \hline
			\end{tabular}
			\caption{Statistics of number of nodes and edges in $G_1\circ G_2.$}
			\label{tab:stat}
		\end{table}

Now we consider counting triads in $G_1\circ G_2.$ Let $s\in\{+, -\}.$ We denote $|E_2^s|^{\stackrel++}$ as the number of edges of sign $s$ which connect two positively marked nodes in $G_2,$ $|E_2^s|^{\pm}$ as the number of edges of sign $s$ which connect one positively marked and one negatively marked nodes in $G_2,$ and $|E_2^s|^{\stackrel--}$ as the number of edges of sign $s$ which connect two negatively marked nodes in $G_2,.$ Then the Table \ref{tab:triad} describes the count of triads of type $T_i$ which denotes a triad having $i$ number of negative edges, $i=0,1,2,3$. 


\begin{table}[H]
			\centering
			\begin{tabular}{l|c|c|c}
				triads & $G_1$ & $G_2$ & $G_1\circ G_2$ \\ \hline\hline
				$\#$ of $T_0$ & $|\T_0(G_1)|$ & $|\T_0(G_2)|$ & $|\T_0(G_1)| + |V_1|\,|\T_0(G_2)| + \M_1^+\,|E_2^+|^{\stackrel++} + \M_1^-|E_2^+|^{\stackrel--}$ \\ 
				$\#$ of $T_1$  & $|\T_1(G_1)|$ & $|\T_1(G_2)|$ & $|\T_1(G_1)| + |V_1|\,|\T_1(G_2)| + \M_1^+\,(|E_2^+|^{\pm}+|E_2^-|^{\stackrel++})$ \\ & & & $+ \M_1^- \, (|E_2^+|^{\pm} + |E_2^-|^{\stackrel--})$  \\ 
				$\#$ of $T_2$ &  $|\T_2(G_1)|$ &  $|\T_2(G_2)|$ & $|\T_2(G_1)| + |V_1|\,|\T_2(G_2)| + \M_1^+\,(|E_2^+|^{\stackrel--} + |E_2^-|^{\pm})$ \\ & & & $+ \M_1^- \, (|E_2^+|^{\stackrel++} + |E_2^-|^{\pm})$ \\ 
$\#$ of $T_3$ & $|\T_3(G_1)|$ & $|\T_3(G_2)|$ & $|\T_3(G_1)| + |V_1|\,|\T_3(G_2)| + \M_1^+\,|E_2^-|^{\stackrel--} + \M_1^-|E_2^-|^{\stackrel++}$ \\ \hline
			\end{tabular}
			\caption{Counts of triads in $G_1\circ G_2.$}
			\label{tab:triad}
		\end{table}

The proof of the formulas presented in Table \ref{tab:stat} and Table \ref{tab:triad} directly follows from the definition of corona product, and easy to verify. Indeed, note that when a node $i$ of $G_1$ gets linked with all the nodes of the $i$th copy of $G_2,$ the total number of new triads created is $|E_2|,$ that is, a triad $(i,j,l)$ is formed for every edge $(j,l)$ in $G_2.$ Now the type of this triad depends on the marking of $i, j, l$ and the sign of the edge $(j,l).$ The signs of all three edges are given by $\mu_1(i)\mu_2(j), \mu_1(i)\mu_2(l)$ and $\sigma_2((j,l)).$ The total number of triads of $G=G_1\circ G_2$ is given by \begin{equation}\label{count:triad}\T(G)= \T(G_1) + |V_1|\,(\T(G_2) + |E_2|)\end{equation} where $\T(G_i)$ denotes the total number of triads in $G_i,$ $i=1,2.$

It is evident that $G_1\circ G_2$ is unbalanced if one of the $G_i,$ $i=1,2$ is unbalanced. However if both $G_i,$ $i=1,2$ are balanced then it is not necessary that $G_1\circ G_2$ is balanced. For instance, the graph $G_1\circ G_2$ in Figure \ref{CP1} is balanced while $G_1$ and $G_2$ are balanced, whereas the graph $G_1\circ G_2$ is unbalanced even if $G_1$ and $G_2$ are balanced in Figure \ref{CP2}. 

\begin{figure}[H]
			\centering			
			\begin{subfigure}{0.2\textwidth}
				\centering
				\begin{tikzpicture}
				\draw [fill] (0, 0) circle [radius=0.1];
				\node [below] at (0, 0) {$1$};
				\draw [fill] (1, 0) circle [radius=0.1];
				\node [below] at (1, 0) {$2$};
				\draw [fill] (2, 0) circle [radius=0.1];
				\node [below] at (2, 0) {$3$};
				\draw (0,0) --(1,0);
				\draw [dashed] (1, 0) -- (2, 0);
				\end{tikzpicture}
				\caption{$G_1$}
			\end{subfigure}
			\begin{subfigure}{0.2\textwidth}
				\centering
				\begin{tikzpicture}
				\draw [fill] (0, 0) circle [radius=0.1];
				\node [below] at (0, 0) {$a$};
				\draw [fill] (1, 0) circle [radius=0.1];
				\node [below] at (1, 0) {$b$};
				\draw [dashed] (0,0) --(1,0);
				\end{tikzpicture}
				\caption{$G_2$}
			\end{subfigure} \\
			\begin{subfigure}{0.2\textwidth}
				\centering
				\begin{tikzpicture}
				\draw [fill] (0, -0.5) circle [radius = .1];
				\node [below] at (0, -0.5) {};
				\node [left] at (0, -0.5) {$b$};
				\draw [fill] (1, 0) circle [radius = .1];
				\node [below] at (1, 0) {$1$};
				\node [above] at (1, 0) {};
				\draw [dashed] (0,-0.5) -- (1, 0);
				\draw [fill] (0, 0.5) circle [radius = .1];
				\node [above] at (0, .5) {};
				\node [left] at (0, .5) {$a$};
				\draw [dashed] (0,-0.5) -- (0, .5);
				\draw [dashed] (0,0.5) -- (1, 0);
				\draw [fill] (2, 1) circle [radius = .1];
				\node [below] at (2, 1) {$2$};
				\node [left] at (2, 1) {};
				\draw (2,1) -- (1, 0);
				\draw [fill] (1.5, 2) circle [radius = .1];
				\node [above] at (1.5, 2) {$a$};
				\node [left] at (1.5, 2) {};
				\draw [dashed] (1.5, 2) -- (2.5, 2);
				\draw [fill] (2.5, 2) circle [radius = .1];
				\node [above] at (2.5, 2) {$b$};
				\node [right] at (2.5, 2) {};
				\draw  (2, 1) -- (1.5, 2);
				\draw  (2, 1) -- (2.5, 2);
				\draw [dashed] (2, 1) -- (3, 0);
				\draw [fill] (3, 0) circle [radius = .1];
				\node [below] at (3, 0) {$3$};
				\node [left] at (3, 0) {};
				\draw  (3, 0) -- (4, -0.5);
				\draw  (3, 0) -- (4, 0.5);
				\draw [dashed] (2, 1) -- (2.5, 2);
				\draw [dashed] (4, 0.5) -- (4, -0.5);
				\draw [fill] (4, -0.5) circle [radius = .1];
				\node [below] at (4, -0.5) {};
				\node [right] at (4,-0.5) {$b$};
				\draw [fill] (4, 0.5) circle [radius = .1];
				\node [above] at (4, 0.5) {};
				\node [right] at (4,0.5) {$a$};
				\end{tikzpicture}
				\caption{$G_1 \circ G_2$ when $G_1, G_2$ are with canonical marking}
			\end{subfigure}\hspace{2.5cm}
	\begin{subfigure}{0.2\textwidth}
				\centering
				\begin{tikzpicture}
				\draw [fill] (0, -0.5) circle [radius = .1];
				\node [below] at (0, -0.5) {};
				\node [left] at (0, -0.5) {$b$};
				\draw [fill] (1, 0) circle [radius = .1];
				\node [below] at (1, 0) {$1$};
				\node [above] at (1, 0) {};
				\draw [dashed] (0,-0.5) -- (1, 0);
				\draw [fill] (0, 0.5) circle [radius = .1];
				\node [above] at (0, .5) {};
				\node [left] at (0, .5) {$a$};
				\draw [dashed] (0,-0.5) -- (0, .5);
				\draw [dashed] (0,0.5) -- (1, 0);
				\draw [fill] (2, 1) circle [radius = .1];
				\node [below] at (2, 1) {$2$};
				\node [left] at (2, 1) {};
				\draw (2,1) -- (1, 0);
				\draw [fill] (1.5, 2) circle [radius = .1];
				\node [above] at (1.5, 2) {$a$};
				\node [left] at (1.5, 2) {};
				\draw [dashed] (1.5, 2) -- (2.5, 2);
				\draw [fill] (2.5, 2) circle [radius = .1];
				\node [above] at (2.5, 2) {$b$};
				\node [right] at (2.5, 2) {};
				\draw [dashed] (2, 1) -- (1.5, 2);
				\draw [dashed] (2, 1) -- (2.5, 2);
				\draw [dashed] (2, 1) -- (3, 0);
				\draw [fill] (3, 0) circle [radius = .1];
				\node [below] at (3, 0) {$3$};
				\node [left] at (3, 0) {};
				\draw  (3, 0) -- (4, -0.5);
				\draw  (3, 0) -- (4, 0.5);
				\draw [dashed] (2, 1) -- (2.5, 2);
				\draw [dashed] (4, 0.5) -- (4, -0.5);
				\draw [fill] (4, -0.5) circle [radius = .1];
				\node [below] at (4, -0.5) {};
				\node [right] at (4,-0.5) {$b$};
				\draw [fill] (4, 0.5) circle [radius = .1];
				\node [above] at (4, 0.5) {};
				\node [right] at (4,0.5) {$a$};
				\end{tikzpicture}
				\caption{$G_1 \circ G_2$ when $G_1, G_2$ are with plurality marking}
			\end{subfigure}
			\caption{The corona product of $G_1\circ G_2$ is shown in (c), (d) with different marking functions on $G_i, i=1,2$.}
			\label{CP2}
		\end{figure}
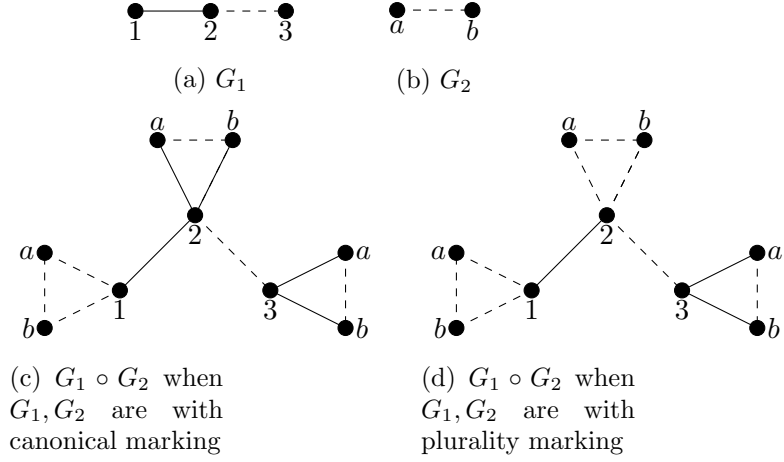	

In the following theorem we classify when $G_1\circ G_2$ will be an unbalanced graph when $G_i$ is balanced, $i=1,2$.
 
\begin{theorem}\label{Thm:b}
Let $G_1=(V_1,E_1,\sigma_1,\mu_1)$ and $G_2=(V_2,E_2,\sigma_2,\mu_2)$ be balanced signed graphs. Then $G_1 \circ G_2$ is unbalanced if and only if $G_2$ contains one of the following types of edges. \begin{itemize}\item[(i)] a positive edge which connects two oppositely marked nodes \item[(ii)] a negative edge which which connects two positively marked nodes \item[(iii)] a negative edge which connects two negatively marked nodes.  \end{itemize}
\end{theorem}
\noindent\pf The proof follows from the fact that any positively or negatively marked node of $G_1$ forms a triad $T_1$ in $G_1\circ G_2$ when there are edges of type $(i)$ and/or $(ii)$ in $G_2,$ otherwise it forms a triad of type $T_3$ when $G_2$ has an edge of type $(iii).$ $\hfill{\square}$

It follows from Theorem \ref{Thm:b} that for two balanced signed graphs  $G_1$ and $G_2,$ the corona product graph $G_1\circ G_2$ is balanced if and only if every positive edge of $G_2$ is incidental to a pair of negatively or positively marked nodes, and every negative edge must be incidental to a pair of oppositely marked nodes.

Now we focus on the spectral properties of $G_1\circ G_2.$ Let $G_1=(V_1,E_1,\sigma_1,\mu_1)$ and $G_2=(V_2,E_2,\sigma_2,\mu_2)$ be two signed graphs with $n$ and $k$ number of nodes, respectively. Suppose $V_1=\{u_1,\hdots,u_{n}\}$ and $V_2=\{v_1,\hdots,v_k\}.$ Let us denote the marking vectors corresponding to vertices in $G_1$ and $G_2$ as $$\mu[V_1]=\bmatrix{\mu_1(u_1) & \mu_1(u_2) & \hdots & \mu_1(u_n)} \,\, \mbox{and} \,\, \mu[V_2]=\bmatrix{\mu_2(v_1) & \mu_2(v_2) & \hdots & \mu_2(v_k)}$$ where $\mu_j(u)=1$ if marking of $u=+,$ otherwise $\mu_j(u)=-1,$ $j=1,2.$

Then with a suitable labeling of the nodes the adjacency matrix of $G_1\circ G_2$ is given by 
\begin{equation} A(G_1\circ G_2)= \bmatrix{A(G_1) & \mu[V_2]\otimes \mbox{diag}(\mu[V_1]) \\ \mu[V_2]^T\otimes \mbox{diag}(\mu[V_1]) & A(G_2)\otimes I_n}\end{equation} where $A(G_i)$ denotes the adjacency matrix associated with $G_i,$ $i=1,2,$ $\otimes$ denotes the Kronecker product of matrices, $I_n$ is the identity matrix of order $n,$ and $$ \mbox{diag}(\mu[V_1]) = \bmatrix{\mu_1(u_1) & 0 & \hdots & 0\\ 0 & \mu_1(u_2) & \hdots & 0\\ \vdots & \vdots & \ddots & \vdots\\ 0 & 0 & \hdots & \mu_1(u_n) }.$$ 

Thus $A(G_1\circ G_2)$ is a symmetric matrix of order $n+nk=n(1+k).$ Recall that a signed graph is called net-regular if $d^+(v) - d^-(v)=d$ is same for every node $v$ of the graph, and $d$ is called the net-regularity of the graph. Then note that the net-regularity of a net-regular graph is always an eigenvalue of the graph and the all-one vector is the corresponding eigenvector. We denote the all-one vector of dimension $k$ as $\textbf{1}_k.$ The following theorem provides adjacency spectra of $G_1\circ G_2$ when $G_2$ is net-regular. 

\begin{theorem}\label{th:adj_spectra_gh}
Let $G_1$ be any signed graph on $n$ nodes and $G_2$ be a net-regular signed graph on $k$ nodes having net-regularity $d.$ Let $(\lam_i, X_i)$ be an adjacency eigenpair of $G_1,$ and $(\eta_j,Y_j)$ be an eigenpair of $G_2,$ $i=1,\hdots,n,$ $j=1,\hdots,k.$ Let $\eta_k=d.$ Then an adjacency eigenpair of $G_1 \circ G_2$ is given by $(\lam_\pm^{(i)}, Z_\pm^{(i)}),$ $i=1,\hdots, n$ where $$\lam_\pm^{(i)} = \frac{d+\lam_i\pm\sqrt{(d-\lam_i)^2 +4k}}{2}, Z_\pm^{(i)} = \bmatrix{X_i \\ \dfrac{\mu_2(v_1)}{\lam_\pm^{(i)}-d}\mbox{diag}(\mu[V_1])X_i \\  \dfrac{\mu_2(v_2)}{\lam_\pm^{(i)}-d}\mbox{diag}(\mu[V_1])X_i \\ \vdots \\  \dfrac{\mu_2(v_k)}{\lam_\pm^{(i)}-d}\mbox{diag}(\mu[V_1])X_i}.$$	

In addition, if all the nodes in $G_2$ are either positively or negatively marked, that is $\mu[V_2]=\textbf{1}_k$ or $-\textbf{1}_k$ then $$\left(\eta_j, \bmatrix{\textbf{0}\\ Y_j\otimes e_i}\right)$$  is an eigenpair of $G_1\circ G_2$ where $j=1,\hdots, k-1,$ and $\{e_i : j=1,\hdots,n\}$ the standard basis of $\R^n.$
 \end{theorem}
\noindent\pf Let $(\lam, Z)$ be an eigenpair of $G_1\circ G_2$ where $\lam\in\R,$ $Z=\bmatrix{Z_0\\ Z_1\\ \vdots \\ Z_k},$ $Z_l\in \R^{n},$ $l=0,1,\hdots, k.$ Then setting $A(G_1\circ G_2)Z=\lam Z$ we obtain $$ \bmatrix{A(G_1) & \mu[V_2]\otimes \mbox{diag}(\mu[V_1]) \\ \mu[V_2]^T\otimes \mbox{diag}(\mu[V_1]) & A(G_2)\otimes I_n}Z=\lam Z$$ which yields the following system of equations. \beano A(G_1)Z_0 +\mu_2(v_1)\mu[V_1]Z_1 + \hdots + \mu_2(v_k)\mu[V_1]Z_k &=& \lam Z_0 \\ \mu_2(v_1)\mu[V_1]Z_0 + \sum_{j=1}^k [A(G_2)]_{1j}Z_j &=& \lam Z_1 \\ &\vdots& \\ \mu_2(v_k)\mu[V_1]Z_0 + \sum_{j=1}^k [A(G_2)]_{2j}Z_j &=& \lam Z_k\eeano where $[A(G_2)]_{xj}$ denote the $(x,j)$th entry of $A(G_2).$ Now adding the last $k$ equations we obtain \beano && (\mu_2(v_1) + \mu_2(v_2) + \hdots +\mu_2(v_k)) \mu[V_1]Z_0 + d(Z_1+\hdots + Z_k) = \lam(Z_1+\hdots+Z_k) \\ &\Rightarrow&  (\mu_2(v_1) + \mu_2(v_2) + \hdots +\mu_2(v_k)) \mu[V_1]Z_0 = (\lam-d)((Z_1+\hdots + Z_k)).\eeano Set $Z_j=\dfrac{\mu_2(v_j)}{\lam-d}\mu[V_1]Z_0$ $j=1,\hdots,k$ and putting these in the first equation of the system we have $$A(G_1)Z_0=\left(\lam - \frac{k}{\lam-d}\right)Z_0.$$ Now setting $Z_0=X_i$ and $\lam - \dfrac{k}{\lam-d}=\lam_i$ we obtain the quadratic polynomial equation  $\lam^2 -\lam(d+\lam_i)-(k-d\lam_i)=0$ solving which the roots $\lam_\pm^{(i)}$ and its corresponding eigenvectors $Z^{(i)}_\pm$ follows. 

If $\mu[V_2]=\textbf{1}_k$ or $-\textbf{1}_k$ then $\mu[V_2]Y_j=0$ since $\textbf{1}_k$ is an eigenvector of $A(G_2),$ and eigenvectors of a symmetric matrix form an orthogonal set. Hence $$\bmatrix{A(G_1) & \mu[V_2]\otimes \mbox{diag}(\mu[V_1]) \\ \mu[V_2]^T\otimes \mbox{diag}(\mu[V_1]) & A(G_2)\otimes I_n}\bmatrix{\textbf{0}\\ Y_j\otimes e_i}=\bmatrix{\mu[V_2]Y_j\otimes \mbox{diag}(\mu[V_1])e_j \\ \eta_j Y_j\otimes e_i }=\eta_j \bmatrix{\textbf{0}\\ Y_j\otimes e_i}.$$ This completes the proof. \hfill{$\square$}
		    	
Note that if $G_2$ is net-regular with plurality marking then either $\mu_2(v_j)=+$ or $\mu_2=-$ for all $j.$ Hence the Theorem \ref{th:adj_spectra_gh} provides all the eigenvalues of $G_1 \circ G_2$ for such a signed graph $G_2.$			

We mention that Barik et al. \cite{barik2007spectrum} gave the complete description of adjacency eigenvalues of $G_1\circ G_2$ when $G_1, G_2$ are unsigned and $G_2$ is regular. Theorem \ref{th:adj_spectra_gh} generalizes their findings, and provides some insights about how signs of the edges and marking of the nodes influence the eigenpairs for corona product graphs.


Now we consider signed Laplacian spectra of corona product graphs. The signed Laplacian matrix associated with a signed graph $G$ on $n$ nodes is defined as $L(G)=D(G)-A(G)$ where $D(G)$ is a diagonal matrix of order $n$ such that the $i$th diagonal entry is $d_i^+ + d_i^-$ and $A(G)$ is the adjacency matrix as usual \cite{zaslavsky1982signed}. The signless Laplacian matrix is defined as  $Q(G)=D(G) + A(G).$ Note that signed Laplacian matrix is a symmetric positive semi-definite matrix, and positive definite when the corresponding graph is unbalanced \cite{kunegis2010spectral,hou2005bounds}. The signless Laplacian matrix associated with a signed graph $G$ and be considered as the signed Laplacian corresponding to the graph $\overline{G}$ which is obtained from $G$ by converting the positive (resp. negative) edges to negative (resp. positive) edges in $G$ \cite{germina2011products}.

First we have the following observation which will be used in sequel. 
\begin{lemma}\label{lem:lev}
Let $G$ be a signed graph on $n$ nodes. Then $\textbf{1}_n$ is an eigenvector corresponding to a signed Laplacian eigenvalue $\lam$ if and only if $d^-_i=d^-(=\lam/2)$ for every node $i$ in $G.$  
\end{lemma}
\noindent\pf From the definition of signed Laplacian matrix, we have $L(G)=D(G)-A(G)$ where $[D(G)]_{ii}=d_i^+ + d^-_i,$ the $i$th diagonal entry of $D(G).$ Then $$L(G)\textbf{1}_n = \bmatrix{2d^-_1\\ 2d^-_2\\ \vdots \\ 2d^-_n}.$$ Hence $\textbf{1}_n$ is an eigenvector of $L(G)$ if and only if $d^-_i$ equals a constant, say $d^-$ for all $i.$ This completes the proof. $\hfill{\square}$

For a pair of signed graphs $G_1=(V_1,E_1,\sigma_1,\mu_1)$ and $G_2=(V_2,E_2,\sigma_2,\mu_2)$ on $n$ and $k$ nodes respectively, the signed and signless Laplacian matrices of $G_1\circ G_2$ are given by \beano L(G_1\circ G_2) &=& \bmatrix{L(G_1)+ kI_n& -\mu[V_2]\otimes \mbox{diag}(\mu[V_1]) \\ -\mu[V_2]^T\otimes \mbox{diag}(\mu[V_1]) & (L(G_2)+I_k)\otimes I_n}\\ Q(G_1\circ G_2) &=& \bmatrix{Q(G_1)+kI_n & \mu[V_2]\otimes \mbox{diag}(\mu[V_1]) \\ \mu[V_2]^T\otimes \mbox{diag}(\mu[V_1]) & (Q(G_2)+I_k)\otimes I_n}\eeano with a suitable labeling on the vertices of $G_1\circ G_2.$ Then we have the following theorem.

\begin{theorem}\label{th:lap1_spectra_gh}
Let $G_1=(V_1,E_1,\sigma_1,\mu_1)$ be a signed graph on $n$ nodes and $G_2=(V_2,E_2,\sigma_2,\mu_2)$ be a signed graph on $k$ nodes. Let $V_2=\{v_1,\hdots,v_k\}.$ Let $(\lam_i, X_i)$ be a signed Laplacian  eigenpair of $G_1,$ $i=1,\hdots,n.$ Let $d^-_j$ denote the negative degree of a node $v_j$ in $G_2.$ Then the roots of the polynomial equations \begin{equation}\label{ev:eqn} x-k-\dfrac{1}{x-(2d_1^- +1)} - \dfrac{1}{x-(2d_2^- +1)} -\hdots - \dfrac{1}{x-(2d_k^- +1)} =\lam_i,\end{equation} $i=1,\hdots,n,$ are signed Laplacian eigenvales of $G_1 \circ G_2.$ An eigenvector corresponding to such an eigenvalue $x=\lam$ corresponding to $\lam_i$ in equation (\ref{ev:eqn})  is given by $$\bmatrix{X_i \\ -\dfrac{\mu_2(v_1)}{\lam-(2d_1^- +1)}\mbox{diag}(\mu[V_1])X_i \\  -\dfrac{\mu_2(v_2)}{\lam-(2d_2^- +1)}\mbox{diag}(\mu[V_1])X_i \\ \vdots \\ -\dfrac{\mu_2(v_k)}{\lam-(2d_k^- +1)}\mbox{diag}(\mu[V_1])X_i}.$$
where $i=1,\hdots, n.$ 
\end{theorem}
\noindent\pf The proof is similar to the proof of Theorem \ref{th:adj_spectra_gh}. Indeed, Let $Z=\bmatrix{Z_0 \\Z_1\\\vdots\\Z_k},$ $Z_j\in\R^n,$ $j=0,\hdots,k$ be an eigenvector corresponding to an eigenvalue $\lam$ of $L(G_1\circ G_2).$ Then from $L(G_1\circ G_2)Z=\lam Z$ we have the following system of equations \beano (L(G_1)+kI_n)Z_0 - \mu_2(v_1)\mbox{diag}(\mu[V_1])Z_1 - \hdots - \mu_2(v_k)\mbox{diag}(\mu[V_1])Z_k &=& \lam Z_0 \\ -\mu_2(v_1)\mbox{diag}(\mu[V_1])Z_0 + ([L(G_2)]_{11} + 1)Z_1 - \hdots + [L(G_2)]_{1k} Z_k &=& \lam Z_1 \\ &\vdots& \\  -\mu_2(v_k)\mbox{diag}(\mu[V_1])Z_0 + [L(G_2)]_{k1} Z_1 - \hdots + ([L(G_2)]_{kk}+1)Z_k &=& \lam Z_k\eeano where $[L(G_2)]_{xy}$ is the $(x,y)$ entry of $L(G_2).$ Then adding the last $k$ equations we have $$-(\mu_2(v_1) + \hdots + \mu_2(v_k))\mbox{diag}(\mu[V_1])Z_0 = \sum_{j=1}^k (\lam -(2d_j^- +1))Z_j.$$  Setting $Z_j=-\dfrac{\mu_2(v_j)}{\lam-(2d^-_j +1)}\mbox{diag}(\mu[V_1]) Z_0,$ $j=1,\hdots,k$ and putting these in the first equation of the system we have $$L(G_1)Z_0 = \left( \lam - k - \dfrac{1}{\lam-(2d_1^- + 1)} - \hdots - \dfrac{1}{\lam-(2d^-_k + 1)} \right)Z_0.$$ Hence the desired result follows by setting $Z_0=X_i.$ $\hfill{\square}$

It is interesting to observe from Theorem \ref{th:adj_spectra_gh} and Theorem\ref{th:lap1_spectra_gh} that adjacency and signed Laplacian eigenvalues of $G_1\circ G_2$ depend only on the degrees of the nodes of $G_2$ whereas the eigenvectors are influenced by degrees of nodes of $G_2$ and the marking scheme of $G_1.$  Besides, observe that the equation (\ref{ev:eqn}) gives the complete list of eigenvalues of $G_1\circ G_2$ only when any two distinct nodes in $G_2$ have distinct negative degrees. Otherwise, the degrees of the polynomial equations (\ref{ev:eqn}) drop and it can not produce all the eigenvalues. In the extreme case when negative degrees of nodes in $G_2$ are all equal then the polynomials become quadratic and produce only $2n$ eigenvalues. The following theorem considers this case.

\begin{theorem}\label{th:lap2_spectra_gh}
Let $G_1=(V_1,E_1,\sigma_1,\mu_1)$ be a signed graph on $n$ nodes and $G_2=(V_2,E_2,\sigma_2,\mu_2)$ be a signed graph on $k$ nodes. Let $V_2=\{v_1,\hdots,v_k\}.$ Let $(\lam_i, X_i)$ be a signed Laplacian  eigenpair of $G_1,$ and $(\eta_j,Y_j)$ are signed Laplacian eigenpairs of $G_2,$ $i=1,\hdots,n,$ $j=1,\hdots,k.$ Let $d^-=d^-_j$ denote the negative degree of every node $v_j$ in $G_2.$ Then a signed Laplacian eigenpair of $G_1 \circ G_2$ is given by $(\lam^{(i)},Z_\pm^{(i)})$ where \beano \lam^{(i)}_\pm &=& \dfrac{2d^-+1+\lam_i+k\pm \sqrt{[(2d^-+1)-(\lam_i+k)]^2+4k }}{2} \\  Z_\pm^{(i)} &=& \bmatrix{X_i \\ -\dfrac{\mu_2(v_1)}{\lam^{(i)}_\pm -(2d^- +1)}\mbox{diag}(\mu[V_1])X_i \\  -\dfrac{\mu_2(v_2)}{\lam^{(i)}_\pm-(2d^- +1)}\mbox{diag}(\mu[V_1])X_i \\ \vdots \\ -\dfrac{\mu_2(v_k)}{\lam^{(i)}_\pm-(2d^- +1)}\mbox{diag}(\mu[V_1])X_i}\eeano where $i=1,\hdots, n.$ 

Let $\eta_k=2d^-.$ In addition if all the nodes in $G_2$ are marked either positively or negatively then an eigenpair of $G_1\circ G_2$ is $$\left(\eta_j + 1, \bmatrix{ 0\\ Y_j\otimes e_i}\right), i=1,\hdots,n, j=1,\hdots,k-1$$ where $\{e_i : i=1,\hdots, n\}$ is the standard basis of $\R^n.$
\end{theorem}
\noindent\pf The first part follows from the proof of Theorem \ref{th:lap1_spectra_gh} by setting $d_j^-=d^-$ for $j=1,\hdots,k.$ Next, if $d_j^-=d^-$ for all $j$ then by Lemma \ref{lem:lev} $\textbf{1}_k$ is an eigenvector of $G_2$ corresponding to the eigenvalue $2d^-=\eta_k.$ Hence $\textbf{1}_k^TY_j=0$ for $j=1,\hdots,k-1$ since $L(G_2)$ is a symmetric matrix and it has complete set of orthogonal vectors. Further, if all the nodes of $G_2$ are positively marked then $\mu[V_2]=\textbf{1}_k^T.$ Finally $$L(G_1\circ G_2)\bmatrix{ 0\\ Y_j\otimes e_i}=\bmatrix{-\mu[V_2]Y_j\otimes e_i \\ (L(G_2)Y_j+Y_j) \otimes e_i}=(\eta_j+1)\bmatrix{0\\ Y_j\otimes e_i}.$$ If the nodes of $G_2$ are marked then $\mu[V_2]=-\textbf{1}_k^T$ and hence the desired result follows. $\hfill{\square}$

We emphasize that the Theorem \ref{th:lap2_spectra_gh} generalizes the result on Laplacian spectra of corona product of unsigned graphs obtained in \cite{barik2007spectrum}. Indeed, setting $d^-=0$ in the Theorem \ref{th:lap2_spectra_gh} describes complete set of eigenpairs of corona product of unsigned graphs. 

\begin{remark} Similar results like Lemma \ref{lem:lev}, Theorem \ref{th:lap1_spectra_gh} and Theorem \ref{th:lap2_spectra_gh} can be obtained for signless Laplacian eigenpairs of $G_1\circ G_2$ only replacing $d_j^-$ by $d_j^+$ for $j=1,\hdots,k.$ Indeed, note that $\textbf{1}_k$ is an eigenvector corresponding to the signless Laplacian eigenvalue $2d^+$ of $G_1\circ G_2$ if and only if $d^+_j=d^+$ for $j=1,\hdots,k.$
\end{remark}

\section{Application to signed network modelling}

Corona product of unsigned graphs has been used as a framework to develop models for complex networks in literature \cite{sharma2017structural,sharma2015self}. Inspired by this idea we propose signed corona graphs as a model for generating large signed networks. In this approach, we consider a small connected signed graph $G=G^{(0)}$ as a seed graph for generation of corona graphs given by $$G^{(m)}=G^{(m-1)}\circ G, m\geq 1$$ with a suitable choice of marking function which will be used in every step $i\leq m-1$ of the formation of $G^{(m)}.$ We explore different structural and spectral properties of $G^{(m)}$ which contribute to investigate how close does it preserve properties of real signed networks. Finally we consider the problem of structural reconstruction of real signed networks using $G^{(m)}$ by a suitable choice of seed graph that can inherit properties of real networks. 

Obviously the structural properties of the seed graph $G$ influence the structural and spectral properties of $G^{(m)}.$ For example, if the edges in $G$ are of the type mentioned in $(i),$ $(ii),$ $(iii)$ of Theorem \ref{Thm:b} and all triads in $G$ are unbalanced then all the triads in $G^{(m)}$ will be unbalanced. Thus in order to have balanced triads in $G^{(m)},$ $G$ must have edges which does not fall in the category of $(i),$ $(ii),$ $(iii)$ of Theorem \ref{Thm:b} or at least one balanced triad in $G$ itself. Below we do a thorough analysis of structural and spectral properties of corona graphs generated by following the canonical marking scheme. We mention that the analytical expression of different properties of corona graphs corresponding to plurality marking is complicated and cumbersome. 

\subsection{Statistics of signed edges and triads in $G^{(m)}$}\label{sec:3.1}

First we have the following observations about $G^{(m)},$ $m\geq 1$ generated by a seed graph $G$ on $n$ nodes with $k$ edges.\begin{enumerate} \item Total number of nodes in $G^{(m)}$ is $n(n+1)^{m}$ \cite{sharma2017structural}. \item The total number of copies of the seed graph in $G^{(m)}$ is $(n+1)^m.$ This follows from that fact that at each step $1\leq i\leq m$ the number of seed graphs added during the formation of $G^{(i)}$ is the number of nodes in $G^{(i-1)}.$ Thus the total number of copies of the seed graph is $$ 1+n+n(n+1)+\hdots+n(n+1)^{m-1} = 1+n\sum_{j=0}^{m-1} (n+1)^j = (n+1)^m. $$
\item The total number of edges added in the formation of $G^{(i+1)}$ from $G^{(i)},$ $0\leq i\leq m-1$ due to the definition of corona product, that is, excluding the edges in $G^{(i)}$ and edges in the newly appeared copies of seed graphs is given by $n^2(n+1)^{i}.$ This follows from that fact that for each node $v$ in $G^{(i)},$ $n$ new edges will be created, and the number of nodes in $G^{(i)}$ is $n(n+1)^{i}.$
\item The total number of edges in $G^{(m)}$ is $k(n+1)^m + \sum_{j=0}^{m-1} n^2(n+1)^j = (k+n)(n+1)^{m} -n.$ \end{enumerate}

We denote $n^{(i)}_+$ and $n^{(i)}_-$ the number of positive and negative nodes in $G^{(i)},$ $0\leq i\leq m-1$ respectively. Thus the total number of nodes in $G^{(i)}=n(n+1)^i (=n^{(i)}_0 + n^{(i)}_- =n^{(i)},$ say). We denote the vertex set which consists of the marked nodes in $G^{(i)}$ as $V^{(i)},$ $ 1\leq i\leq m$. Let the number of positive and negative links in $G^{(0)}$ be denoted as $e_+^{(0)}$ and $e_-^{(0)}$ respectively.

Consider counting the number of signed edges in $G^{(m)}.$ Note that the number of copies of the seed graph added during the formation of $G^{(i+1)}$ from $G^{(i)}$ is $n(n+1)^i$ (which is exactly the number of nodes in $G^{(i)}$), $0\leq i\leq m-1.$ Then in $G^{(i+1)}$ the number of new positive and negative links are created as \begin{eqnarray} e^{(i+1)}_+ &=& e_+^{(0)}n(n+1)^i + n^{(i)}_+ \, n_+^{(0)} + n_-^{(i)} \, n_-^{(0)} \label{ne:eqn1}\\ e_-^{(i+1)} &=& e_-^{(0)}n(n+1)^i + n_+^{(i)} \, n_-^{(0)} + n_-^{(i)} \, n_+^{(0)}\label{ne:eqn2}\end{eqnarray} respectively. The first term in the expression of $e^{(i+1)}_+$ counts the number of positive links which are part of the newly added copies of seed graph, the second term is the number of links which are formed by joining positively marked nodes in $V^{(i)}$ and positively marked nodes in the copies of the seed graph, and the third term is the number of positive edges formed by existing negatively marked nodes in $G^{(i)}$ and the newly appeared negatively marked nodes. Similarly, the expression of $e^{(i+1)}_-$ follows. 

Now let us count the number of positively and negatively marked nodes in $G^{(i+1)}$ after its formation from $G^{(i)}.$ There can be two types of nodes in $G^{(i+1)},$ the newly appeared nodes $v\in V^{(i+1)}\setminus V^{(i)}$ and the existing set of nodes $V^{(i)}$ of $G^{(i)}.$ First, consider the nodes in $V^{(i+1)} \setminus V^{(i)}.$ The number of positively and negatively marked nodes among such nodes shall be $n_+^{(i)}n$ and $n^{(i)}_- n$ respectively. This follows from the fact that when a new node (positively/ negatively marked) gets linked to an existing positively (resp. negatively) marked node in $V^{(i)}$ then it becomes positively (resp. negatively) marked due to the definition of canonical marking. Now let us consider the marking of nodes in $V^{(i+1)}\cap V^{(i)}.$ Note that any node $v\in V^{(i+1)}\cap V^{(i)}$ was already marked during the formation of $G^{(i)}.$ When $G^{(i+1)}$ is formed, each such node forms links to all the nodes in a copy of a newly appeared seed graph. If it was $\mu(v)=+$ then it remains the same marking when $n^{(0)}_-$ is even and $n^{(0)}_+$ is even or odd, whereas it changes its marking only when $n_-^{(0)}$ is odd and $n^{(0)}_+$ is either even or odd. Similarly, if it was $\mu(v)=-$ then it remains the same marking for even value of $n_+^{(0)}$ and any value of $n^{(0)}_-,$ whereas it changes its marking when $n_+^{(0)}$ is odd and any value of $n^{(0)}_-.$ Thus finally we have the following. The number of total positively and negatively marked nodes in $G^{(i+1)}$ is described in Table \ref{tab:nmstat}.

\begin{table}[H]
\centering
\begin{tabular}{l|c|c|c|c}
Marking & Condition on $G$ & in $V^{(i+1)}\setminus V^{(i)}$ &  in $V^{(i+1)}\cap V^{(i)}$ &  in $G^{(i+1)}$ \\ \hline\hline
\multirow{4}{*}{$+$} & $n^{(0)}_+$ odd  & \multirow{4}{*}{$n^{(i)}_+ \,n$}  & \multirow{2}{*}{$n^{(i)}_-$} & \multirow{2}{*}{$n^{(i)}_- + n\, n^{(i)}_+ $}\\ 
&  $n^{(0)}_-$ odd & &  &\\
\cline{2-2}\cline{4-5}
& $n^{(0)}_+$ even  &  & \multirow{2}{*}{$0$} & \multirow{2}{*}{$n^{(i)}_+ \,n$}\\ 
& $n^{(0)}_-$ odd & & &\\
\cline{2-2}\cline{4-5}
& $n^{(0)}_+$ even  &  & \multirow{2}{*}{$n_+^{(i)}$} & \multirow{2}{*}{$n_+^{(i)} (1 + n)$}\\ 
& $n^{(0)}_-$ even & & &\\
\cline{2-2}\cline{4-5}
& $n^{(0)}_+$ odd  &  & \multirow{2}{*}{$n^{(i)}$} & \multirow{2}{*}{$n^{(i)} +n^{(i)}_+ n$}\\ 
& $n^{(0)}_-$ even & & &\\
\hline\hline
\multirow{4}{*}{$-$} & $n^{(0)}_+$ odd  & \multirow{4}{*}{$n^{(i)}_- \,n$}  & \multirow{2}{*}{$n^{(i)}_+$} & \multirow{2}{*}{$n^{(i)}_+ +  n\,n^{(i)}_- $}\\ 
&  $n^{(0)}_-$ odd & &  &\\
\cline{2-2}\cline{4-5}
& $n^{(0)}_+$ even  &  & \multirow{2}{*}{$n^{(i)}$} & \multirow{2}{*}{$n^{(i)}_- \,n + n^{(i)}$}\\ 
& $n^{(0)}_-$ odd & & &\\
\cline{2-2}\cline{4-5}
& $n^{(0)}_+$ even  &  & \multirow{2}{*}{$n_-^{(i)}$} & \multirow{2}{*}{$n_-^{(i)} (1 + n)$}\\ 
& $n^{(0)}_-$ even & & &\\
\cline{2-2}\cline{4-5}
& $n^{(0)}_+$ odd  &  & \multirow{2}{*}{$0$} & \multirow{2}{*}{$n^{(i)}_- \, n$}\\ 
& $n^{(0)}_-$ even & & &\\
\hline \hline
\end{tabular}
\caption{The relation between the number of marked nodes in $G^{(i+1)}$ and $G^{(i)}.$  The last column displays the values of $n^{(i+1)}_+$ and $n^{(i+1)}_-.$}
\label{tab:nmstat}
\end{table}

Now we are in a position to estimate the number of signed links in $G^{(m)}$ by utilizing equations (\ref{ne:eqn1}),  (\ref{ne:eqn2}) and the fact that $n_+^{(i)}=n^{(i)} - n_-^{(i)}.$
 \beano e^{(m)}_+ &=& e^{(0)}_++ \sum_{i=0}^{m-1} e^{(i+1)}_+ =e^{(0)}_+\, \left(1 + n \sum_{i=0}^{m-1} (n+1)^i\right) + n^{(0)}_+ \sum_{i=0}^{m-1} n_+^{(i)} + n_-^{(0)} \sum_{i=0}^{m-1} n_-^{(i)}\\ &=& e^{(0)}_+ (n+1)^{m} + n^{(0)}_+[(n+1)^{m} -1] + (n_-^{(0)} - n_+^{(0)}) \sum_{i=0}^{m-1} n_-^{(i)}\eeano where the values of $n_-^{(i)}$ can be computed from the recurrence relation given in Table \ref{tab:nmstat} for all $i.$ Similarly, \beano e^{(m)}_- &=& e^{(0)}_- (n+1)^{m} + n^{(0)}_-[(n+1)^{m} -1] + (n_+^{(0)} - n_-^{(0)}) \sum_{i=0}^{m-1} n_-^{(i)}.\eeano 

It is evident from the above expressions that marking of nodes in the seed graph plays the key role. Besides, the number of positive and negative edges become equal when the number of positively and negatively marked nodes are same in the seed graph.  Now we focus on couting the number of triads of type $T_i,$ $i=0,1,2,3$ in $G^{(m)}.$ First observe that there can be six types of links in a signed graph. A positive/negative link can be incidental to positively marked nodes, oppositely marked nodes or negatively marked nodes. When such nodes in a triad in a copy of the seed graph get linked during the formation of each step of the corona graph $G^{(i+1)}$ from $G^{(i)},$ differnt types of triads are created. The Figure \ref{triads} depicts all the possible cases.

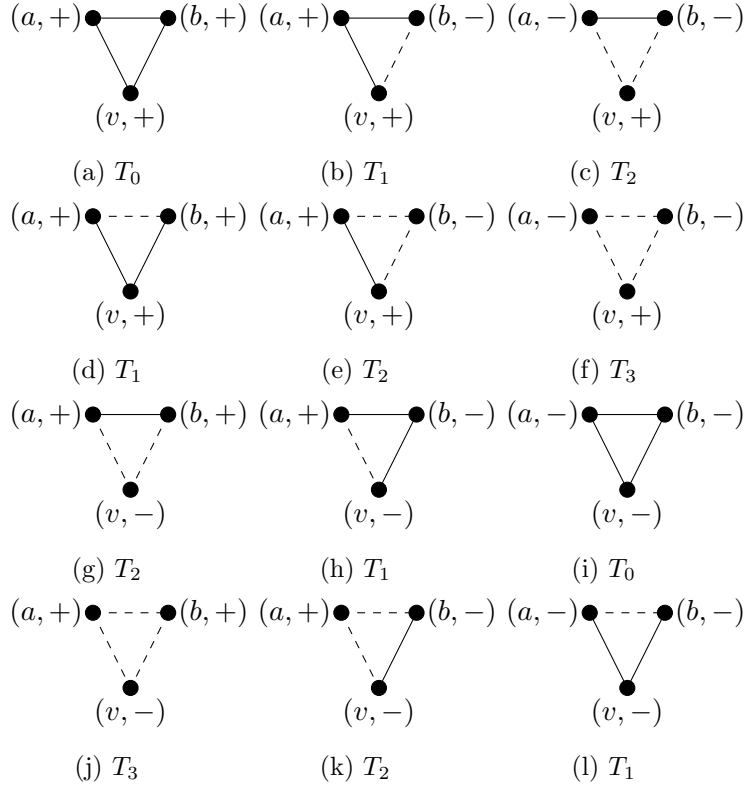
\begin{figure}[H]
			\centering			
			\begin{subfigure}{0.2\textwidth}
				\centering
				\begin{tikzpicture}
				\draw [fill] (0, 0) circle [radius=0.1];
				\node [left] at (0, 0) {$(a,+)$};
				\draw [fill] (1, 0) circle [radius=0.1];
				\node [right] at (1, 0) {$(b,+)$};
				\draw (0,0) --(1,0);
                                          \node [below] at (0.5, -1) {$(v,+)$};
				\draw [fill] (0.5, -1) circle [radius=0.1];
                                           \draw (0,0) --(0.5,-1);
				\draw (1, 0) -- (0.5, -1);
				\end{tikzpicture}
				\caption{$T_0$}
			\end{subfigure}
                              \hspace{0.2cm}
                              \begin{subfigure}{0.2\textwidth}
				\centering
				\begin{tikzpicture}
				\draw [fill] (0, 0) circle [radius=0.1];
				\node [left] at (0, 0) {$(a,+)$};
				\draw [fill] (1, 0) circle [radius=0.1];
				\node [right] at (1, 0) {$(b,-)$};
				\draw (0,0) --(1,0);
                                          \node [below] at (0.5, -1) {$(v,+)$};
				\draw [fill] (0.5, -1) circle [radius=0.1];
                                           \draw (0,0) --(0.5,-1);
				\draw [dashed](1, 0) -- (0.5, -1);
				\end{tikzpicture}
				\caption{$T_1$}
			\end{subfigure}
                              \hspace{0.2cm}
                              \begin{subfigure}{0.2\textwidth}
				\centering
				\begin{tikzpicture}
				\draw [fill] (0, 0) circle [radius=0.1];
				\node [left] at (0, 0) {$(a,-)$};
				\draw [fill] (1, 0) circle [radius=0.1];
				\node [right] at (1, 0) {$(b,-)$};
				\draw (0,0) --(1,0);
                                          \node [below] at (0.5, -1) {$(v,+)$};
				\draw [fill] (0.5, -1) circle [radius=0.1];
                                           \draw [dashed](0,0) --(0.5,-1);
				\draw [dashed](1, 0) -- (0.5, -1);
				\end{tikzpicture}
				\caption{$T_2$}
			\end{subfigure} \\
			\begin{subfigure}{0.2\textwidth}
				\centering
				\begin{tikzpicture}
				\draw [fill] (0, 0) circle [radius=0.1];
				\node [left] at (0, 0) {$(a,+)$};
				\draw [fill] (1, 0) circle [radius=0.1];
				\node [right] at (1, 0) {$(b,+)$};
				\draw [dashed](0,0) --(1,0);
                                          \node [below] at (0.5, -1) {$(v,+)$};
				\draw [fill] (0.5, -1) circle [radius=0.1];
                                           \draw (0,0) --(0.5,-1);
				\draw (1, 0) -- (0.5, -1);
				\end{tikzpicture}
				\caption{$T_1$}
			\end{subfigure}
                              \hspace{0.2cm}
                              \begin{subfigure}{0.2\textwidth}
				\centering
				\begin{tikzpicture}
				\draw [fill] (0, 0) circle [radius=0.1];
				\node [left] at (0, 0) {$(a,+)$};
				\draw [fill] (1, 0) circle [radius=0.1];
				\node [right] at (1, 0) {$(b,-)$};
				\draw [dashed](0,0) --(1,0);
                                          \node [below] at (0.5, -1) {$(v,+)$};
				\draw [fill] (0.5, -1) circle [radius=0.1];
                                           \draw (0,0) --(0.5,-1);
				\draw [dashed](1, 0) -- (0.5, -1);
				\end{tikzpicture}
				\caption{$T_2$}
			\end{subfigure}
                              \hspace{0.2cm}
                              \begin{subfigure}{0.2\textwidth}
				\centering
				\begin{tikzpicture}
				\draw [fill] (0, 0) circle [radius=0.1];
				\node [left] at (0, 0) {$(a,-)$};
				\draw [fill] (1, 0) circle [radius=0.1];
				\node [right] at (1, 0) {$(b,-)$};
				\draw [dashed](0,0) --(1,0);
                                          \node [below] at (0.5, -1) {$(v,+)$};
				\draw [fill] (0.5, -1) circle [radius=0.1];
                                           \draw [dashed](0,0) --(0.5,-1);
				\draw [dashed](1, 0) -- (0.5, -1);
				\end{tikzpicture}
				\caption{$T_3$}
			\end{subfigure} \\
			\begin{subfigure}{0.2\textwidth}
				\centering
				\begin{tikzpicture}
				\draw [fill] (0, 0) circle [radius=0.1];
				\node [left] at (0, 0) {$(a,+)$};
				\draw [fill] (1, 0) circle [radius=0.1];
				\node [right] at (1, 0) {$(b,+)$};
				\draw (0,0) --(1,0);
                                          \node [below] at (0.5, -1) {$(v,-)$};
				\draw [fill] (0.5, -1) circle [radius=0.1];
                                           \draw [dashed](0,0) --(0.5,-1);
				\draw [dashed](1, 0) -- (0.5, -1);
				\end{tikzpicture}
				\caption{$T_2$}
			\end{subfigure}
                              \hspace{0.2cm}
                              \begin{subfigure}{0.2\textwidth}
				\centering
				\begin{tikzpicture}
				\draw [fill] (0, 0) circle [radius=0.1];
				\node [left] at (0, 0) {$(a,+)$};
				\draw [fill] (1, 0) circle [radius=0.1];
				\node [right] at (1, 0) {$(b,-)$};
				\draw (0,0) --(1,0);
                                          \node [below] at (0.5, -1) {$(v,-)$};
				\draw [fill] (0.5, -1) circle [radius=0.1];
                                           \draw [dashed](0,0) --(0.5,-1);
				\draw (1, 0) -- (0.5, -1);
				\end{tikzpicture}
				\caption{$T_1$}
			\end{subfigure}
                              \hspace{0.2cm}
                              \begin{subfigure}{0.2\textwidth}
				\centering
				\begin{tikzpicture}
				\draw [fill] (0, 0) circle [radius=0.1];
				\node [left] at (0, 0) {$(a,-)$};
				\draw [fill] (1, 0) circle [radius=0.1];
				\node [right] at (1, 0) {$(b,-)$};
				\draw (0,0) --(1,0);
                                          \node [below] at (0.5, -1) {$(v,-)$};
				\draw [fill] (0.5, -1) circle [radius=0.1];
                                           \draw (0,0) --(0.5,-1);
				\draw (1, 0) -- (0.5, -1);
				\end{tikzpicture}
				\caption{$T_0$}
			\end{subfigure} \\
			\begin{subfigure}{0.2\textwidth}
				\centering
				\begin{tikzpicture}
				\draw [fill] (0, 0) circle [radius=0.1];
				\node [left] at (0, 0) {$(a,+)$};
				\draw [fill] (1, 0) circle [radius=0.1];
				\node [right] at (1, 0) {$(b,+)$};
				\draw [dashed](0,0) --(1,0);
                                          \node [below] at (0.5, -1) {$(v,-)$};
				\draw [fill] (0.5, -1) circle [radius=0.1];
                                           \draw [dashed](0,0) --(0.5,-1);
				\draw [dashed](1, 0) -- (0.5, -1);
				\end{tikzpicture}
				\caption{$T_3$}
			\end{subfigure}
                              \hspace{0.2cm}
                              \begin{subfigure}{0.2\textwidth}
				\centering
				\begin{tikzpicture}
				\draw [fill] (0, 0) circle [radius=0.1];
				\node [left] at (0, 0) {$(a,+)$};
				\draw [fill] (1, 0) circle [radius=0.1];
				\node [right] at (1, 0) {$(b,-)$};
				\draw [dashed](0,0) --(1,0);
                                          \node [below] at (0.5, -1) {$(v,-)$};
				\draw [fill] (0.5, -1) circle [radius=0.1];
                                           \draw [dashed](0,0) --(0.5,-1);
				\draw (1, 0) -- (0.5, -1);
				\end{tikzpicture}
				\caption{$T_2$}
			\end{subfigure}
                              \hspace{0.2cm}
                              \begin{subfigure}{0.2\textwidth}
				\centering
				\begin{tikzpicture}
				\draw [fill] (0, 0) circle [radius=0.1];
				\node [left] at (0, 0) {$(a,-)$};
				\draw [fill] (1, 0) circle [radius=0.1];
				\node [right] at (1, 0) {$(b,-)$};
				\draw [dashed](0,0) --(1,0);
                                          \node [below] at (0.5, -1) {$(v,-)$};
				\draw [fill] (0.5, -1) circle [radius=0.1];
                                           \draw (0,0) --(0.5,-1);
				\draw (1, 0) -- (0.5, -1);
				\end{tikzpicture}
				\caption{$T_1$}
			\end{subfigure}
\caption{Formation of signed triads in every step of the corona product. Each signed edge consisiting nodes $(a,p), (b,q)$ joins an exiting node $(v,r)$ where $p,q,r\in\{+,-\}$ denote the marking of the nodes.}
			\label{triads}
		\end{figure}

We denote $|E^{(0)}_s|^{\stackrel pq}$ as the number of links in $G^{(0)}$ having sign $s$ that is adjacent to nodes with markings $p$ and $q.$ For example, $|E^{(0)}_-|^{\stackrel +-}$ denotes the number of negative links in $G^{(0)}$ adjacent to oppositely marked nodes. Then the Table \ref{tab:triad} gives the number of triads of $T_j, j=0,1,2,3$ added after the formation of $G^{(i+1)}$ from $G^{(i)}$ due to the new signed edges created between the nodes in $V^{(i)}$ and $V^{(i+1)}\setminus V^{(i)},$ $1\leq i\leq m-1.$ We denote this number as $\T_j^{(i+1)}.$

\begin{table}[H]
\centering
\begin{tabular}{c|c}
$j$ & $\T_j^{(i+1)}$ \\
\hline
$j=0$ & $n_+^{(i)} \, |E^{(0)}_+|^{\stackrel ++} + n_-^{(i)} \, |E^{(0)}_+|^{\stackrel --}$ \\
$j=1$ & $n_+^{(i)} \,\left( |E^{(0)}_+|^{\stackrel +-} +  |E^{(0)}_-|^{\stackrel ++} \right) + n_-^{(i)} \,\left( |E^{(0)}_+|^{\stackrel +-} +  |E^{(0)}_-|^{\stackrel --}\right)$ \\
$j=2$ & $n_+^{(i)} \, \left( |E^{(0)}_+|^{\stackrel --} +  \, |E^{(0)}_-|^{\stackrel +-} \right) + n_-^{(i)} \, \left( |E^{(0)}_+|^{\stackrel ++} +  \, |E^{(0)}_-|^{\stackrel +-}\right)$ \\
$j=3$ & $n_+^{(i)} \, |E^{(0)}_-|^{\stackrel --} + n_-^{(i)} \, |E^{(0)}_-|^{\stackrel ++}$ \\
\hline
\end{tabular}
\caption{The number of $T_j$ type triads between the nodes in $V^{(i)}$ and $V^{(i+1)}\setminus V^{(i)}$ created during the formation of $G^{i+1}$ from $G^{(i)}.$}
\label{tab:gmstat}
\end{table}

Note that if any type of triad is present in a seed graph $G=G^{(0)},$ it will appear every time a copy of the seed graph gets attached to the formation of $G^{(m)}.$ Let $\T^{(0)}_j$ denote the number of triads of type $T_j$ that exist in $G^{(0)}.$  Then the total number of $T_j$ $j=0,1,2,3$ type triads $\T_j(G^{(m)})$ in $G^{(m)}$ can be calculated as follows.
\beano \T_0(G^{(m)}) &=& \T^{(0)}_0\, \left(1 + n \sum_{i=0}^{m-1} (n+1)^i\right) + \sum_{i=0}^{m-1} \T_0^{(i+1)} \\ &=& \T^{(0)}_0 (n+1)^{m} + |E^{(0)}_+|^{\stackrel ++} \sum_{i=0}^{m-1} n_+^{(i)}  +  |E^{(0)}_+|^{\stackrel --} \sum_{i=0}^{m-1} n_-^{(i)}  \\ &=& \T^{(0)}_0 (n+1)^{m} + |E^{(0)}_+|^{\stackrel ++} [(n+1)^m -1] + \left(|E^{(0)}_+|^{\stackrel --} - |E^{(0)}_+|^{\stackrel ++}\right) \sum_{i=0}^{m-1} n_-^{(i)}\eeano where the first term includes the number of $T_0$ type triads which are in the copies of the seed graph  in $G^{(m)},$ and the remaining terms are due to the new triads which are born during the formation of the corona product in each step of the formation of $G^{(m)}.$ The values of $n_-^{(i)}$ can be found from Table \ref{tab:nmstat} for each $i.$ Similarly we have the following. 

\beano \T_1(G^{(m)}) &=&  \T^{(0)}_1 (n+1)^{m} + \left( |E^{(0)}_+|^{\stackrel +-} +  |E^{(0)}_-|^{\stackrel ++} \right)  [(n+1)^m -1] + \left(  |E^{(0)}_-|^{\stackrel --} -  |E^{(0)}_-|^{\stackrel ++} \right) \sum_{i=0}^{m-1} n_-^{(i)} \\
\T_2(G^{(m)}) &=&   \T^{(0)}_2 (n+1)^{m} +\left( |E^{(0)}_+|^{\stackrel --} +  \, |E^{(0)}_-|^{\stackrel +-} \right) [(n+1)^m -1]  +  \left( |E^{(0)}_+|^{\stackrel ++} -  \, |E^{(0)}_+|^{\stackrel --} \right) \sum_{i=0}^{m-1} n_-^{(i)} \\
\T_3(G^{(m)}) &=&   \T^{(0)}_3 (n+1)^{m} + |E^{(0)}_-|^{\stackrel --}  [(n+1)^m -1] + \left(  |E^{(0)}_-|^{\stackrel ++} -  |E^{(0)}_-|^{\stackrel --} \right) \sum_{i=0}^{m-1} n_-^{(i)}\eeano where the values of $n_-^{(i)}$ can be computed by using the recurrence relation given in Table \ref{tab:nmstat}.

Note from the above formulas that the number of $T_j$ type triads in $G^{(m)}$ can be controlled by choosing the seed graph with an appropriate distribution of signed edges whose adjacent nodes are marked in a particular fashion. 

\subsection{Degrees of nodes of $G^{(m)}$}

Note that after the appearance of a node in a particular step during the formation of $G^{(m)}$, from next step onward the total degree of the node increases by $n$ at every step until the process stops. The dynamics of the change in positive and negative degrees depend on the existing marking of the node and the marking of the nodes which get linked to it at every step. Since the new nodes which get linked to an existing node are the nodes in the copy of the seed graph, the change in degrees of the existing node depends on the marking of the nodes in the seed graph. We first derive the positive and negative degrees of the nodes in $G=G^{(0)},$ and then we consider nodes in $G^{(i)}$ that is which appear at the $i$th step during the formation of $G^{(m)},$ $1\leq i\leq m.$ The initial marking of a node $v$ is denoted by $\mu(v).$ By initial marking we mean once the node $v$ joins the existing graph.

We have the following theorem which provides positive and negative degree of a node in the seed graph $G^{(0)}$ after the formation of $G^{(m)}.$

\begin{theorem}
Let $G^{(m)}, m\geq 1$ be the signed corona graph corresponding to a seed graph $G=G^{(0)}$ on $n>1$ nodes defined by the canonical marking $\mu.$ Let the positive and negative degree of $v$ in $G^{(0)}$ be denoted by $d^+_0(v)$ and $d^-_0(v)$ respectively. Suppose $G^{(0)}$ has $n_+$ positively marked nodes and $n_-$ negatively marked nodes. Then for any $v\in V^{(0)}$ the positive degree $d^+(v)$ and negative degree $d^-(v)$ of $v$ in $G^{(m)}$ are given by the Table \ref{tab:deg0}.
\begin{table}[H]
\centering
\begin{tabular}{l|c|c|c|c}
 Condition  & $\mu(v)$  & $d^+(v)$  &  $d^-(v)$ & Condition \\ 
on $G$ &  &   &   & on $m$ \\ \hline\hline
$n_+$ odd, &  \multirow{4}{*}{$+$}  & \multirow{2}{*}{$d^+_0(v)+\frac{m}{2}(n_++n_-)$}  & \multirow{2}{*}{$d^-_0(v)+\frac{m}{2}(n_+ +n_-)$} & \multirow{2}{*}{even}\\ 
$n_-$ odd &&&& \\
 &  & \multirow{2}{*}{$d^+_0(v)+\frac{m+1}{2}n_++ \frac{m-1}{2}n_-$} & \multirow{2}{*}{$d^-_0(v)+\frac{m+1}{2}n_-+\frac{m-1}{2}n_+$} & \multirow{2}{*}{odd}\\ 
&&&&\\
\cline{2-5}
&  \multirow{4}{*}{$-$}  &  \multirow{2}{*}{$d^+_0(v)+\frac{m}{2}(n_++n_-)$}  &  \multirow{2}{*}{$d^-_0(v)+\frac{m}{2}(n_+ +n_-)$} &  \multirow{2}{*}{even}\\
&&&&\\ 
&  &  \multirow{2}{*}{$d^+_0(v)+\frac{m+1}{2}n_-+ \frac{m-1}{2}n_+$} &  \multirow{2}{*}{$d^-_0(v)+\frac{m+1}{2}n_++\frac{m-1}{2}n_-$} & odd\\
&&&&\\ 
\hline\hline
$n_+$ even, &   \multirow{2}{*}{$+$}  &  \multirow{2}{*}{$d^+_0(v)+n_++(m-1)n_-$}  &  \multirow{2}{*}{$d^-_0(v)+n_-+(m-1)n_+$} &  \multirow{2}{*}{even/odd}\\
$n_-$ odd &&&&\\ 
\cline{2-5}
  &  \multirow{2}{*}{$-$}  & \multirow{2}{*}{$d^+_0(v)+m n_-$}  & \multirow{2}{*}{$d^-_0(v)+m n_+$} & \multirow{2}{*}{even/odd}\\
&&&&\\ 
\hline \hline
$n_+$ even, &  \multirow{2}{*}{$+$}  & \multirow{2}{*}{$d^+_0(v) + mn_+$}  & \multirow{2}{*}{$d^-_0(v) + m n_-$} & \multirow{2}{*}{even/odd}\\
$n_-$ even &&&&\\ 
\cline{2-5}
  &   \multirow{2}{*}{$-$}  &  \multirow{2}{*}{$d^+_0(v)+m n_-$}  &  \multirow{2}{*}{$d^-_0(v)+m n_+$} &  \multirow{2}{*}{even/odd}\\ 
&&&&\\
\hline  \hline
$n_+$ odd, &   \multirow{2}{*}{$+$}  &  \multirow{2}{*}{$d^+_0(v) + mn_+$}  &  \multirow{2}{*}{$d^-_0(v) + m n_-$} &  \multirow{2}{*}{even/odd}\\
$n_-$ even &&&&\\ 
\cline{2-5}
  &   \multirow{2}{*}{$-$}  &  \multirow{2}{*}{$d^+_0(v) + n_- + (m-1) n_+$}  &  \multirow{2}{*}{$d^-_0(v) + n_+ + (m-1) n_-$} &  \multirow{2}{*}{even/odd}\\ 
&&&&\\
\hline
\end{tabular}
\caption{Positive and negative degree of nodes of the seed graph $G$ in $G^{(m)}$}
\label{tab:deg0}
\end{table}
\end{theorem}
\noindent\pf Let $n_+$ and $n_-$ be odd. If a node $v$ in $G$ is positively marked then at every step of the formation of $G^{(m)}$ the marking of the node will change since at every step it will generate odd number of negative links to the newly appeared nodes in a copy of the seed graph. If $\mu(v)=+$ (resp. $\mu(v)=-$) at any step then at the next step it will generate $n_-$ (resp. $n_+$) negative links. If $m$ is even then at $m/2$ number of steps it will be incidental to $n_-$ and $n_+$ nodes. Thus the desired result follows. A similar argument can prove the case for $m$ being odd.  If $\mu(v)=-$ it's marking shall change at every step and hence the desired result follows. 

Let $n+$ be even and $n_-$ odd. If $\mu(v)=+$ then next step the $v$ will generate $n_-$ negative links and since $n_-$ is odd, its marking become $-.$ From next step onward it will generate even number $n_+$ of negative links at every step as $n_+$ is even. Hence $\mu(v)$ remains $-$ for all the steps till $m.$ Thus the desired follows. If $\mu(v)=-,$ it remains negatively marked for all the steps since it generates even number of negative links. This completes the proof. 

Let $n_+$ and $n_-$ both be even. If $\mu(v)=+$, it remains positively marked for all the steps since every time it generates $n_-$ number of negative links and $n_-$ is even. Similarly, if $\mu(v)=-,$ it remains negatively marked as every times it creates even number of new negative links. Hence the proof. 

Finally let $n_+$ be odd and $n_-$ even. If $\mu(v)=+$ it remains $+$ since at every step it creates even number $n_-$ of negative links. If $\mu(v)=-$ then at the next step the marking becomes $+$ since it creates $n_+$ number of negative links and $n_+$ is odd. Next step onward it creates even number $n_-$ of negative links and hence the marking of $v$ remains $+$ for ever until $m$th step. Hence the desired result follows. $\hfill{\square}$

Now we consider the positive and negative degrees of a node $v$ which appears at the $i$th step, that is, during the formation of  $G^{(i)}$, $1\leq i\leq m$ and gets linked with an existing node $u$ in $G^{(i-1)}.$ Obviously $u$ is a node in a copy of the seed graph and it has a marking $\mu(v)$ as a node in the seed graph but as soon as it joins the node $u$ in $G^{(i-1)},$ the marking of $v$ may change depending on the marking of the node $u.$ Let us denote the degree of $v$ as $d^\pm_0(v)$ when we consider it as a node of the seed graph and the positive or negative degree increases by $1$ when it joins $u$ in $G^{(i-1)}$ depending on the marking $\mu(u).$ The following theorem provides the positive and negative degree of such a node $v$.

\begin{theorem}
Let $v$ be a node which joins at the $i$th step, $1\leq i\leq m$ of the formation of the corona graph $G^{(m)}$ generated by a seed graph $G=G^{(0)}$ which contains $n_+$ positively marked and $n_-$ negatively marked nodes. Let $v$ gets attached with a node $u$ in $G^{(i-1)}$ during its birth in $G^{(i)}.$ Then the positive degree $d^+(v)$ and negative degree of $d^-(v)$ in $G^{(m)}$ are given by the Table \ref{tab:degi}. We denote the marking of $v$ as $\mu_0(v)$ when we consider it as a node in the copy of a seed graph which joins the graph $G^{(i-1)},$ whereas $\mu(u)$ denotes the marking of $u$ as a node in $G^{(i-1)}.$ 

\begin{table}[H]
\centering 
\begin{tiny}
\begin{tabular}{l|c|c|c|c|c}
 Condition  & $(\mu(u),\mu_0(v))$ & $d^+(v)$  &  $d^-(v)$ & $m$ & $i$ \\ 
on $G$ &  &  &     & & \\ \hline\hline
$n_+$ odd, & $(+,+)$  &  \multirow{2}{*}{$d^+_0(v)+1+\frac{m-i}{2}n$}  &  \multirow{2}{*}{$d^-_0(v)+\frac{m-i}{2}n$} & odd & odd\\ 
$n_-$ odd & &  & & even & even \\ \cline{3-6}
 & &   \multirow{2}{*}{$d^+_0(v)+1+\frac{m-i+1}{2}n_++ \frac{m-i-1}{2}n_-$} & \multirow{2}{*}{$d^-_0(v)+\frac{m-i+1}{2}n_-+\frac{m-i-1}{2}n_+$} & odd & even \\ 
 & &  &   & even & odd \\ \cline{2-6}
 &  $(+,-)$&  \multirow{2}{*}{$d^+_0(v)+\frac{m-i}{2}n$}  &  \multirow{2}{*}{$d^-_0(v) + 1+\frac{m-i}{2}n$} & odd & odd\\ 
 & & &  & even & even \\ \cline{3-6}
 & &   \multirow{2}{*}{$d^+_0(v)+\frac{m-i+1}{2}n_++ \frac{m-i-1}{2}n_-$} & \multirow{2}{*}{$d^-_0(v)+1+\frac{m-i+1}{2}n_-+\frac{m-i-1}{2}n_+$} & odd & even \\ 
 & &  &   & even & odd \\ \cline{2-6}
 &  $(-,+)$&  \multirow{2}{*}{$d^+_0(v)+\frac{m-i}{2}n$}  &  \multirow{2}{*}{$d^-_0(v) + 1+\frac{m-i}{2}n$} & odd & odd\\ 
 & & &  & even & even \\ \cline{3-6}
 & &   \multirow{2}{*}{$d^+_0(v)+\frac{m-i+1}{2}n_-+ \frac{m-i-1}{2}n_+$} & \multirow{2}{*}{$d^-_0(v)+1+\frac{m-i+1}{2}n_++\frac{m-i-1}{2}n_-$} & odd & even \\ 
 & &  &   & even & odd \\ \cline{2-6}
 &  $(-,-)$&  \multirow{2}{*}{$d^+_0(v)+1+\frac{m-i}{2}n$}  &  \multirow{2}{*}{$d^-_0(v) +\frac{m-i}{2}n$} & odd & odd\\ 
 & & &  & even & even \\ \cline{3-6}
 & &   \multirow{2}{*}{$d^+_0(v)+1+\frac{m-i+1}{2}n_-+ \frac{m-i-1}{2}n_+$} & \multirow{2}{*}{$d^-_0(v)+\frac{m-i+1}{2}n_++\frac{m-i-1}{2}n_-$} & odd & even \\ 
 & &  &   & even & odd \\ 
\hline \hline
$n_+$ even, &  \multirow{2}{*}{$(+,+)$}  &   \multirow{2}{*}{$d^+_0(v)+1+n_++(m-i-1)n_-$ if $m>i$}   &  \multirow{2}{*}{$d^-_0(v)+n_-+(m-i-1)n_+$ if $m>i$} &  \multirow{2}{*}{any} &  \multirow{2}{*}{any}\\ 
 & &  \multirow{2}{*}{$d^+_0(v)+1$ if $m=i$} &  \multirow{2}{*}{$d^-_0(v)$ if $m=i$}& &  \\ 
$n_-$ odd & &  & & &  \\ \cline{2-6}
 &  \multirow{2}{*}{$(+,-)$}  &   \multirow{2}{*}{$d^+_0(v)+n_++(m-i-1)n_-$ if $m>i$}  &    \multirow{2}{*}{$d^-_0(v)+1+n_-+(m-i-1)n_+$ if $m>i$} &  \multirow{2}{*}{any} &  \multirow{2}{*}{any}\\ 
 & &  \multirow{2}{*}{$d^+_0(v)$ if $m=i$} &  \multirow{2}{*}{$d^-_0(v)+1$ if $m=i$}& &  \\ 
 & &  & & &  \\ \cline{2-6}
 &  \multirow{2}{*}{$(-,+)$}  &  \multirow{2}{*}{$d^+_0(v)+(m-i)n_-$}   &   \multirow{2}{*}{$d^-_0(v)+1+(m-i)n_+$} &  \multirow{2}{*}{any} &  \multirow{2}{*}{any}\\ 
 & &  & & &  \\ \cline{2-6}
 &  \multirow{2}{*}{$(-,-)$}  &  \multirow{2}{*}{$d^+_0(v)+1+(m-i)n_-$}  &  \multirow{2}{*}{$d^-_0(v)+(m-i)n_+$} &  \multirow{2}{*}{any} &  \multirow{2}{*}{any}\\ 
 & &  & & &  \\ 
\hline\hline
$n_+$ even, &  \multirow{2}{*}{$(+,+)$}  &  \multirow{2}{*}{$d^+_0(v)+1+(m-i)n_+$}  &   \multirow{2}{*}{$d^-_0(v)+(m-i)n_-$} &  \multirow{2}{*}{any} &  \multirow{2}{*}{any}\\ 
$n_-$ even & &  & & &  \\ \cline{2-6}
 &  \multirow{2}{*}{$(+,-)$}  &  \multirow{2}{*}{$d^+_0(v)+(m-i)n_+$}  &   \multirow{2}{*}{$d^-_0(v)+1+(m-i)n_-$} &  \multirow{2}{*}{any} &  \multirow{2}{*}{any}\\ 
 & &  & & &  \\ \cline{2-6}
 &  \multirow{2}{*}{$(-,+)$}  &  \multirow{2}{*}{$d^+_0(v)+(m-i)n_-$}  &   \multirow{2}{*}{$d^-_0(v)+1+(m-i)n_+$} &  \multirow{2}{*}{any} &  \multirow{2}{*}{any}\\ 
 & &  & & &  \\ \cline{2-6}
 &  \multirow{2}{*}{$(-,-)$}  &  \multirow{2}{*}{$d^+_0(v)+1+(m-i)n_-$}  &   \multirow{2}{*}{$d^-_0(v)+(m-i)n_+$} &  \multirow{2}{*}{any} &  \multirow{2}{*}{any}\\ 
 & &  & & &  \\ 
\hline\hline
$n_+$ odd, &  \multirow{2}{*}{$(+,+)$}  &  \multirow{2}{*}{$d^+_0(v)+1+(m-i)n_+$}  &   \multirow{2}{*}{$d^-_0(v)+(m-i)n_-$} &  \multirow{2}{*}{any} &  \multirow{2}{*}{any}\\ 
$n_-$ even & &  & & &  \\ \cline{2-6}
 &  \multirow{2}{*}{$(+,-)$}  &  \multirow{2}{*}{$d^+_0(v)+(m-i)n_+$}  &   \multirow{2}{*}{$d^-_0(v)+1+(m-i)n_-$} &  \multirow{2}{*}{any} &  \multirow{2}{*}{any}\\ 
 & &  & & &  \\ \cline{2-6}
 &  \multirow{2}{*}{$(-,+)$}  &  \multirow{2}{*}{$d^+_0(v)+n_-+(m-i-1)n_+$ if $m> i$}  &   \multirow{2}{*}{$d^-_0(v)+1+n_++(m-i-1)n_-$  if $m> i$} &  \multirow{2}{*}{any} &  \multirow{2}{*}{any}\\ 
 & & \multirow{2}{*}{$d^+_0(v)$ if $m=i$}   & \multirow{2}{*}{$d^-_0(v)+1$ if $m=i$} & & \\
&&&&& \\ \cline{2-6}
 &  \multirow{2}{*}{$(-,-)$}  &   \multirow{2}{*}{$d^+_0(v)+1+n_-+(m-i-1)n_+$ if $m> i$} &    \multirow{2}{*}{$d^-_0(v)+n_++(m-i-1)n_-$  if $m> i$}  &  \multirow{2}{*}{any} &  \multirow{2}{*}{any}\\ 
 & & \multirow{2}{*}{$d^+_0(v)+1$ if $m=i$}   & \multirow{2}{*}{$d^-_0(v)$ if $m=i$} & & \\
 & &  & & &  \\ 
\hline
\end{tabular}
\end{tiny}
\caption{Positive and negative degree of a node which appears at the $i$th step during the formation of $G^{(m)}, 1\leq i\leq m$ }
\label{tab:degi}
\end{table}
\end{theorem}
\noindent\pf The proof follows by observing the change of marking of a node $v$ due to the appearance of $n$ new nodes which connect to $v$ in every step after it appears in the $i$th step of the formation of $G^{(m)}.$ Let $v$ be a node which appears in $G^{(i)}$, $1\leq i\leq m.$ Recall that a node is positively marked if and only if it is attached to even number of negatively signed links, otherwise the marking of a node is negative. 

First assume that $n_-$ and $n_+$ are odd. Let $\mu(u)=+$ and $\mu_0(v)=+.$ Note that sign of the edge between $u$ and $v$ is $\mu(u) \times \mu_0(v)=+.$ Hence the marking of $v$ in $G^{(i)}$ is $+$ as the number of negative edges adjacent to it does not increase by $1$ and $\mu_0(v)=+$ implies it was attached to even number of negative edges. At the $(i+1)$th step $v$ forms $n_-$ negative edges and $n_-$ is odd. Hence $\mu(v)=-$ in $G^{(i+1)}.$ Proceeding similarly, the sequence of marking of $v$ will be $+, -, +, -, \hdots.$ Thus the positive degree of $v$ increases as follows. At $G^{(i)}$ the degree of $v$ becomes $d_0^+(v) +1.$ At $(i+1)$th step the positive degree of $v$ becomes $d^+_0(v) + 1 + n_+,$ and the negative degree becomes $d_0^-(v)+n_-.$ Since marking of $v$ is $-$ at $(i+1)$th step, it will create $n_-$ positive edges and $n_+$ negative edges at the next step. Thus the positive and negative degree will be $(d_0^+(v) +1)+n_+ +n_-+n_+ + \hdots$ and the negative degree will be $d_0^-(v)+n_-+n_++n_-+\hdots,$ and both the series stop after $m-i$ steps. The desired result follows by considering different cases of $m$ and $i$ being  even and/ odd. 

If $\mu(u)=+$ and $\mu_0(v)=-$ the pattern of change of markings of $v$ starting from the $i$ th step is $+, -, +, -, \hdots.$ The increase of positive and negative degree are given by $d^+_0(v)+n_+ +n_-+n_++n_-+\hdots$ and $(d^-_0+1)+n_-+n_++n_-+n_++\hdots.$ Thus the result follows by considering different cases for $m$ and $i.$

Similarly it can be easily verified that the positive and negative degree of $v$ will be $d^+_0(v)+n_-+n_++n_-+\hdots$ and $(d^-_0(v) +1) +n_++n_-+n_++\hdots$ respectively when $\mu(u)=-$ and $\mu_0(v)=+.$ Finally if $\mu(u)=-$ and $\mu_0(v)=-$ then the positive and negative degree of $v$ will be $(d^+_0(v)+1)+n_-+n_++n_-+\hdots$ and $d_0^- + n_++n_-+n_++\hdots$ respectively. Thus the desired result follows.

Following a similar procedure for different values of $n_-$ and $n_+$ all the desired result can be verified. This completes the proof. $\hfill{\square}$

It is evident from the expressions of the positive and negative degree from Table \ref{tab:deg0} and Table \ref{tab:degi} that as $i$ increases the degree decreases. This means a degree of node is higher than the degree of another node if it joins the network before than the another one. Finally the degrees of the nodes in $G^{(0)}$ are the highest ones.  

For the computation of degree distribution of corona graphs the frequency or total number of occurrence of a node with a particular degree need to be determined. Thus the question is how many nodes do have a particular degree given in the tables above? It is obvious from the construction of the corona graph that the number of nodes which appear at step $i$ with a given degree $d^\pm(v)$ with $(\mu(u),\mu_0(v))=(+,+)$ is $n^{(i)}_+n_+;$ $(\mu(u),\mu_0(v))=(+,-)$ is $n^{(i)}_+n_-;$ $(\mu(u),\mu_0(v))=(-,+)$ is $n^{(i)}_-n_+;$ and $(\mu(u),\mu_0(v))=(-,-)$ is $n^{(i)}_-n_-,$ for $1\leq i\leq m.$ Below we compute positive and negative degree distributions of corona graphs $G^{(i)}, 0\leq i\leq 5$ in Figure \ref{fig:pn_deg_dist} for the seed graphs $G_1, G_2$.   

\begin{figure}[H]
			\begin{subfigure}[b]{0.5\textwidth}
				\centering
				\begin{tikzpicture}
				\draw [fill] (0, 0) circle [radius=0.1];
				\draw [fill] (1, 0) circle [radius=0.1];
				\draw [fill] (0.5, 1) circle [radius=0.1];
				\draw [dashed] (0,0) -- (0.5,1);
				\draw (1,0) -- (0.5,1);
				\end{tikzpicture}
				\caption{$G_1$}
			\end{subfigure}				
			\begin{subfigure}[b]{0.5\textwidth}
				\centering
				\begin{tikzpicture}
				\draw [fill] (0, 0) circle [radius=0.1];
				\draw [fill] (1, 0) circle [radius=0.1];
				\draw [fill] (0,1) circle [radius=0.1];
				\draw [fill] (1,1) circle [radius=0.1];
				\draw (1,0) --(1,1);
				\draw [dashed] (0,0)	-- (1,0);
				\draw (0,0) -- (0,1);
				\draw (0,1) -- (1,1);
				\end{tikzpicture}
				\caption{$G_2$}
			\end{subfigure}	
\\
			\begin{subfigure}[b]{0.5\textwidth}
				\centering
				\includegraphics[scale=0.4]{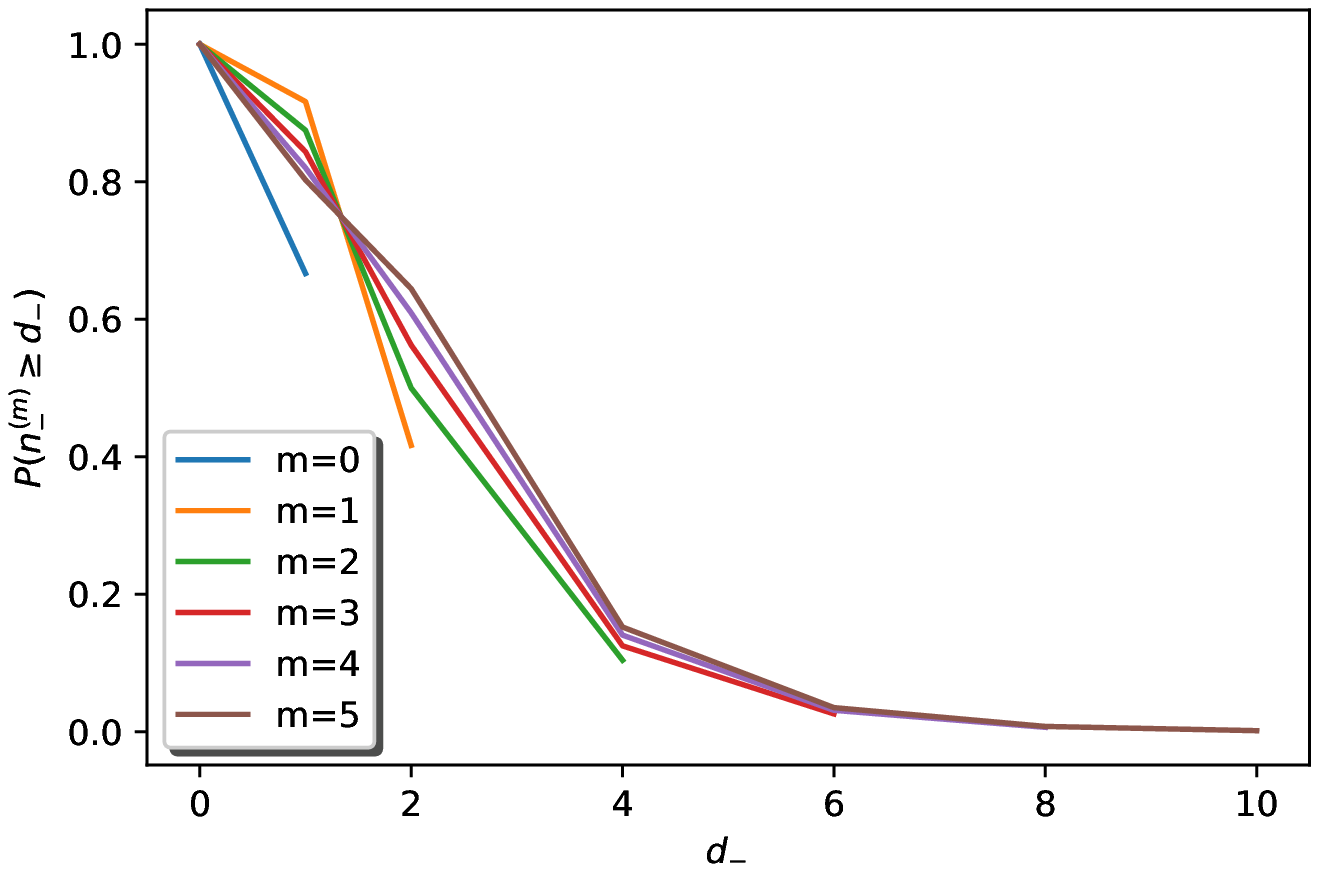}
				\caption{Negative degree distribution of $G_1^{(m)}$}
			\end{subfigure}
                              \begin{subfigure}[b]{0.5\linewidth}
				\centering
				\includegraphics[scale=0.4]{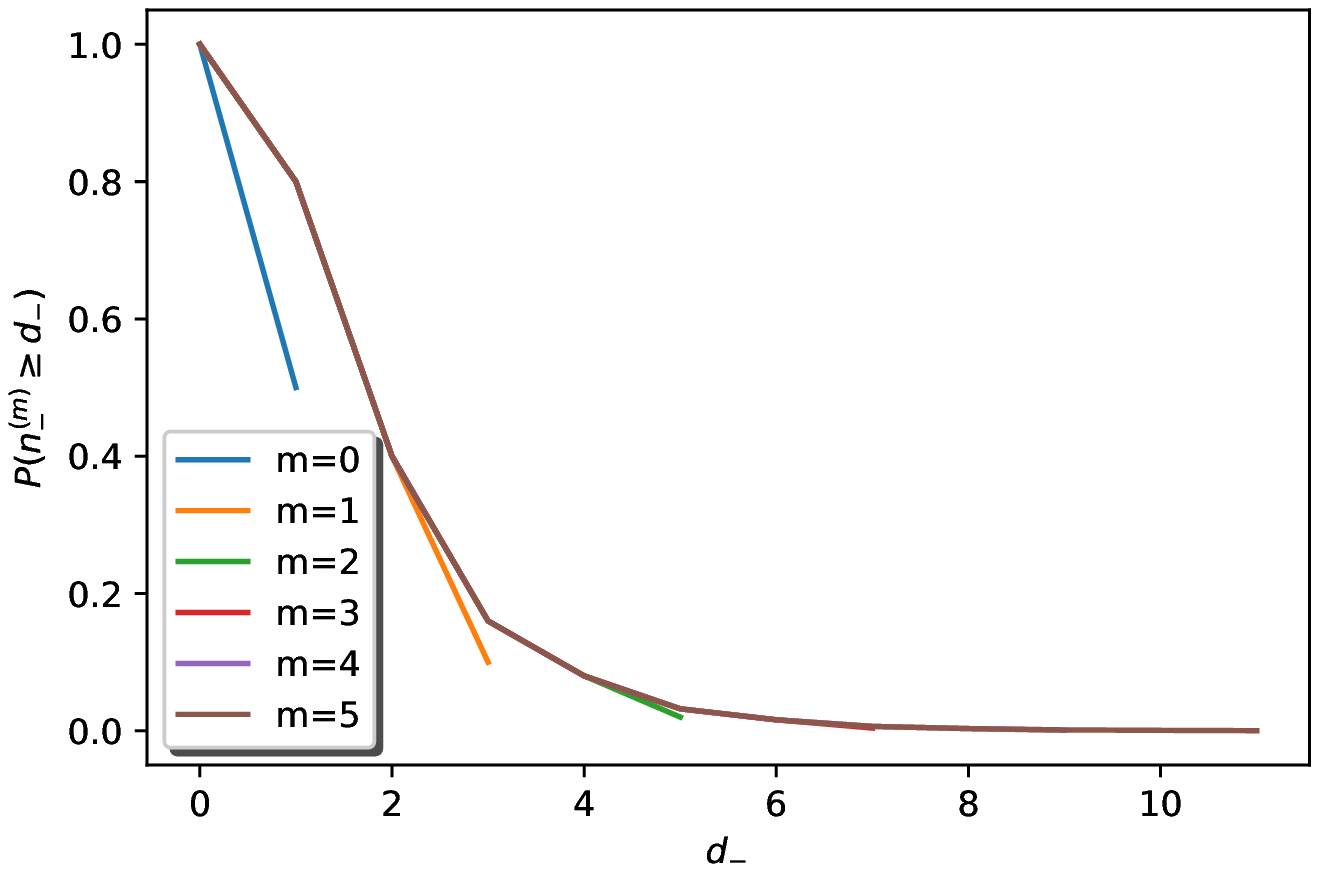}
				\caption{Negative degree distribution of $G_2^{(m)}$}
			\end{subfigure}
			\begin{subfigure}[b]{0.5\linewidth}
				\centering
				\includegraphics[scale=0.4]{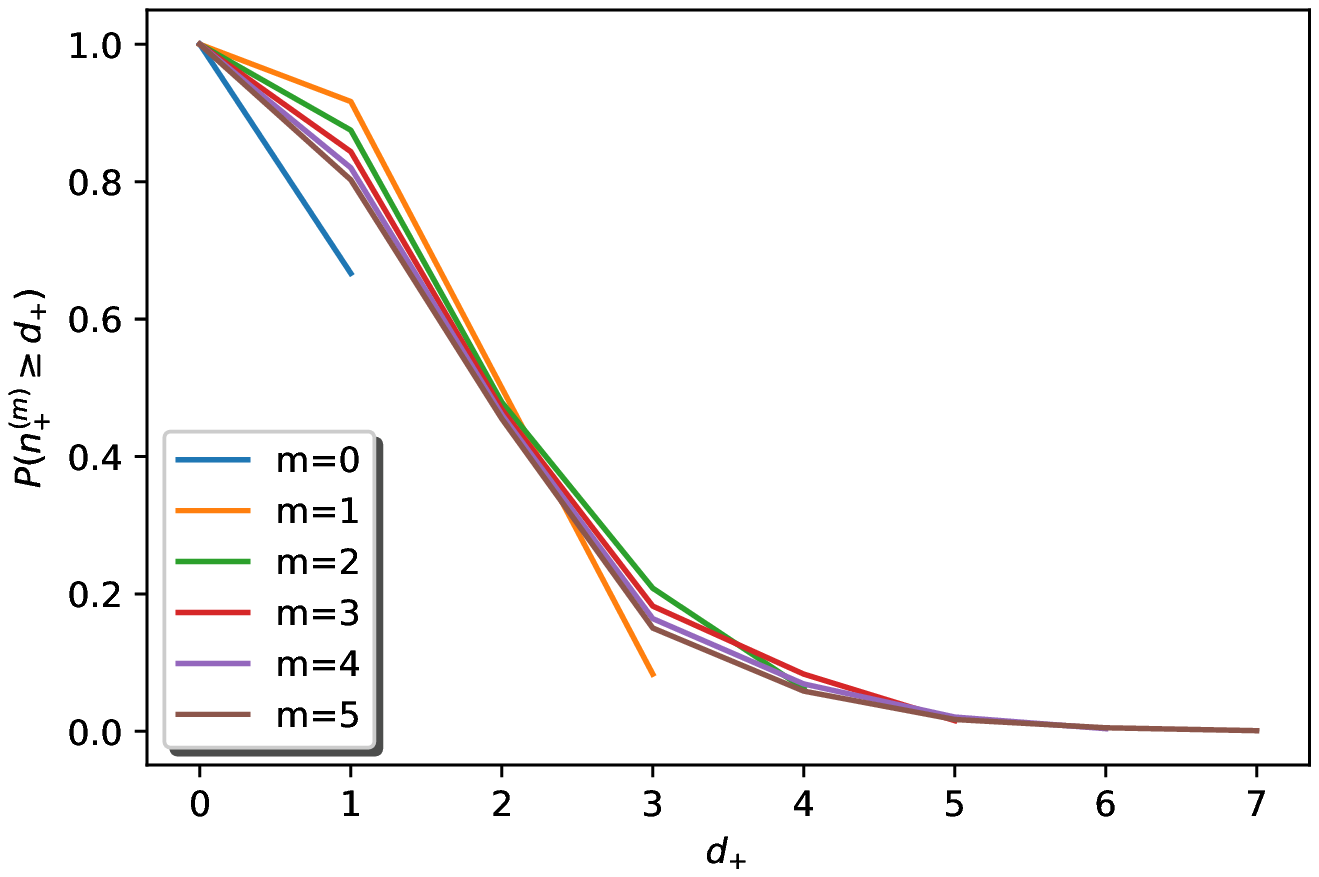}
				\caption{Positive degree distribution of $G_1^{(m)}$}
			\end{subfigure}
			\begin{subfigure}[b]{0.5\linewidth}
				\centering
				\includegraphics[scale=0.4]{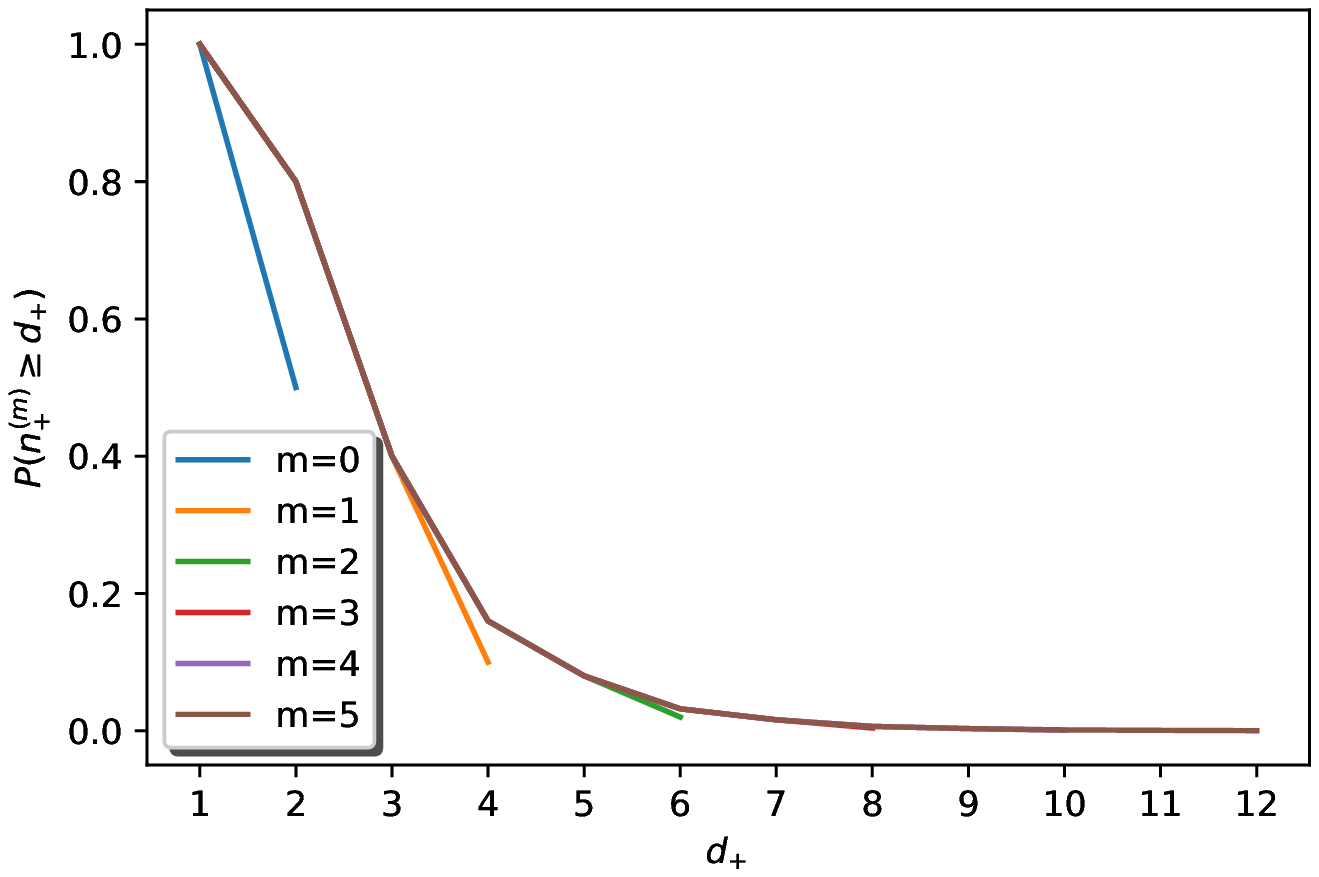}
				\caption{Positive degree distribution of $G_2^{(m)}$}
			\end{subfigure}
			\caption{Degree distributions of $G_i^{(m)}$, $i=1,2$.} 
			\label{fig:pn_deg_dist}
		\end{figure}

\subsection{Algebraic conflict of corona graphs}

Recall that algebraic conflict of a signed graph is the least signed Laplacian eigenvalue of the graph \cite{aref2019signed,kunegis2010spectral}. Theorem \ref{th:lap1_spectra_gh} gives the complete list of signed Laplacian eigenvalues only when distinct nodes of the corresponding graph have distinct negative degrees. However for a large graph the computation of these eigenvalues by calculating roots of high degree polynomials is computationally challenging. Due to the exponential growth of corona graphs, withing a few steps the corona graphs contain huge number of nodes even if the seed graph has considerably small number of nodes. Thus we consider Theorem  \ref{th:lap2_spectra_gh} which gives the collection of all signed Laplacian eigenvalues of the signed graph $G_1\circ G_2$ when marking of all the nodes of $G_2$ are either $+$ or $-,$ and negative degree of all the nodes in $G_2$ are same. Hence setting a seed graph $G=G^{(0)}$ having these properties we provide expression of signed Laplacian eigenvalues of $G^{(m)}$ as follows by iterative use of Theorem \ref{th:lap2_spectra_gh}.

First we have the following lemma which will be used in sequel.

\begin{lemma}\label{lem:mn}
Let $n,m$ be positive integers. Then $$2^mn+ \sum_{i=1}^{m-1}2^in(n-1)(n+1)^{m-i-1} + n(n-1)(n+1)^{m-1}=n(n+1)^m.$$
\end{lemma}
\pf The proof follows from Lemma 1 in \cite{sharma2017structural}.\hfill{$\square$}

Consider the function $f:\R\rightarrow \R$ defined by \begin{equation}\label{def:f}f(x)=\frac{2d^-+1+x+n\pm\sqrt{[(2d^-+1)-(x+n)]^2+4n}}{2}\end{equation} where $d^-$ is a nonnegative integer and $n$ is a positive integer. Denote $f^i$ as $i$ number of compositions of $f$ with itself. For example, $f^1=f$ and $f^2=f\circ f.$ Then we have the following theorem.

\begin{theorem}\label{Thm:gmev}
Let $G$ be a signed graph on $n$ nodes such that $d^-$ is the fixed negative degree of each node of $G$ and every node of $G$ is either positively or negatively marked. Suppose $\eta_j, j=1,\hdots,n$ are the signed Laplacian eigenvalues of $G.$ Let $\eta_n=2d^-.$ Then signed Laplacian eigenvalues of $G^{(m)}, m\geq 1$ defined by the seed graph $G$ are given by \begin{itemize} \item[(a)] $f^{m}(\eta_j)$ with multiplicity $1,$ $j=1,\hdots, n$ \item[(b)] $f^i(\eta_j+1)$ with multiplicity $n(n+1)^{m-i-1},$ $1\leq i\leq m-1,$ $j=1,\hdots,n-1,$ and $m>1$\item[(c)] $\eta_j+1$ with multiplicity $n(n+1)^{m-1}, j=1,\hdots,n-1.$ \end{itemize} 
\end{theorem}
\pf First note that $2d^-$ is a signed Laplacian eigenvalue of $G$ due to Lemma \ref{lem:lev}. The expression of signed laplacian eigenvalues of $G^{(m)}$ follows from iterative use of Theorem \ref{th:lap2_spectra_gh}. Using Lemma \ref{lem:mn} it can be easily verified that this is the complete list of signed Laplacian eigenvalues of $G^{(m)}.$ \hfill{$\square$}

The algebraic conflict of corona graphs defined by a seed graph satisfying the properties of Theorem \ref{Thm:gmev} can be calculated by taking the minimum of the signed Laplacian eigenvalues of $G^{(m)}$ as given in Theorem \ref{Thm:gmev}. However it is interesting to verify that $2d^-$ is a fixed point of the function $$f(x)=\frac{2d^-+1+x+n - \sqrt{[(2d^-+1)-(x+n)]^2+4n}}{2}$$ for any $n.$ This implies that if $2d^-$ is the algebraic conflict of a seed graph $G$ with every node having negative degree $d^-,$ $2d^-$ remains the algebraic conflict of $G^{(m)}$ for any $m\geq 1.$ An example of a seed graph having algebraic conflict $2d^-$ is given by 
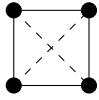
\begin{figure}[H]
			\centering	
				\begin{tikzpicture}
				\draw [fill] (0, 0) circle [radius=0.1];
				\draw [fill] (1, 0) circle [radius=0.1];
				\draw [fill] (0,1) circle [radius=0.1];
				\draw [fill] (1,1) circle [radius=0.1];
				\draw (1,0) --(1,1);
				\draw (0,0)	-- (1,0);
				\draw (0,0) -- (0,1);
				\draw (0,1) -- (1,1);
				\draw [dashed] (0,1) -- (1,0);
				\draw [dashed] (0,0) -- (1,1);
				\end{tikzpicture}
				\caption{Seed graph $G$ with algebraic conflict $2$}
		\end{figure}
where $n=4$ and $d^-=1.$

Otherwise, for seed graphs $G^{(0)}$ if $2d^-$ is not a signed Laplacian eigenvalue , we conjecture that $\eta_{\min}+1 \leq f^i(\eta_j+1), 1\leq i\leq m-1, 1\leq j\leq n-1$ and $\eta_{\min}+1\leq f^m(\eta_j), 1\leq j\leq n,$ and hence algebraic conflict of $G^{(m)}$ is algebraic conflict of $G^{(0)}$ plus one. We verify the same for all possible signed graphs on $3,4,5$ nodes satisfying the conditions in Theorem \ref{Thm:gmev}. In that case it would be a remarkable feature that degree of balance of corona graphs $G^{(m)}$ is not affected by values of $m$ and large signed corona graphs can have small algebraic conflict which is at the order of the algebraic conflict of the seed graph. We consider a few seed graphs on $3,4,5$ nodes that are given in Figure \ref{fig:d} and compute the algebraic conflict of $G^{(m)}, 1\leq m\leq 5$ in Figure \ref{fig:ac}. 


\begin{figure}[H]
			\centering	
\begin{subfigure}[b]{0.15\textwidth}
				\centering
				\begin{tikzpicture}
				\draw [fill] (0, 0) circle [radius=0.1];
				\draw [fill] (1, 0) circle [radius=0.1];
				\draw [fill] (0.5,1) circle [radius=0.1];
				\draw [dashed] (0,0) -- (1,0);
				\draw [dashed] (0,0) -- (0.5,1);
				\draw [dashed] (1,0) -- (0.5,1);
				\end{tikzpicture}
				\caption{$G_1$}
			\end{subfigure}
					\begin{subfigure}[b]{0.15\textwidth}
				\centering
				\begin{tikzpicture}
				\draw [fill] (0, 0) circle [radius=0.1];
				\draw [fill] (1, 0) circle [radius=0.1];
				\draw [fill] (0,1) circle [radius=0.1];
				\draw [fill] (1,1) circle [radius=0.1];
				\draw [dashed] (1,0) --(1,1);
				\draw [dashed] (0,0)	-- (1,0);
				\draw [dashed] (0,0) -- (0,1);
				\draw [dashed] (0,1) -- (1,1);
				\draw (0,1) -- (1,0);
				\draw (0,0) -- (1,1);
				\end{tikzpicture}
				\caption{$G_2$}
			\end{subfigure}		
\begin{subfigure}[b]{0.15\textwidth}
				\centering
				\begin{tikzpicture}
				\draw [fill] (0, 0) circle [radius=0.1];
				\draw [fill] (1, 0) circle [radius=0.1];
				\draw [fill] (0,1) circle [radius=0.1];
				\draw [fill] (1,1) circle [radius=0.1];
                                           \draw [fill] (2,0.5) circle [radius=0.1];
				\draw (1,0) --(1,1);
				\draw [dashed] (0,0)	-- (1,0);
				\draw [dashed] (0,0) -- (0,1);
                                          \draw [dashed] (1,0) -- (2,0.5);
                                          \draw [dashed] (1,1) -- (2,0.5);
                                          \draw (0,1) -- (2,0.5);
				\draw [dashed] (0,1) -- (1,1);
				\draw (0,0) -- (2,0.5);
				\end{tikzpicture}
				\caption{$G_3$}
			\end{subfigure}
\hspace{0.2cm}
\begin{subfigure}[b]{0.15\textwidth}
				\centering
				\begin{tikzpicture}
				\draw [fill] (0, 0) circle [radius=0.1];
				\draw [fill] (1, 0) circle [radius=0.1];
				\draw [fill] (0,1) circle [radius=0.1];
				\draw [fill] (1,1) circle [radius=0.1];
                                           \draw [fill] (2,0.5) circle [radius=0.1];
				\draw (1,0) --(1,1);
				\draw [dashed] (0,0)	-- (1,0);
				\draw [dashed] (0,0) -- (0,1);
                                          \draw [dashed] (1,0) -- (2,0.5);
                                          \draw [dashed] (1,1) -- (2,0.5);
                                          \draw (0,0) -- (2,0.5);
                                          \draw (0,1) -- (2,0.5);
				\draw [dashed] (0,1) -- (1,1);
				\draw (0,1) -- (1,0);
				\draw (0,0) -- (1,1);
				\end{tikzpicture}
				\caption{$G_4$}
			\end{subfigure}
			\caption{Sample seed graphs on $3,4,5$ nodes.}
			\label{fig:d}
		\end{figure}
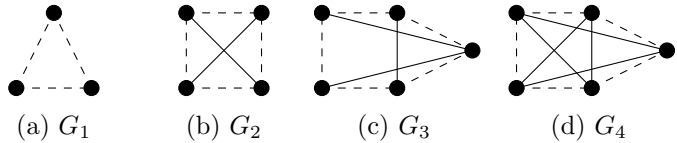

\begin{figure}[H]
				\centering
				\includegraphics[scale=0.4]{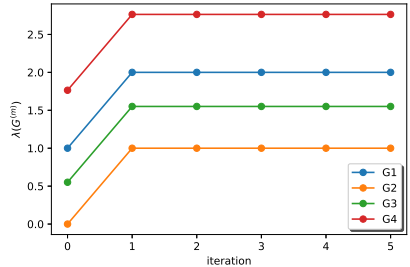}
			\caption{Algebraic conflicts of corona graphs generated by seed graphs listed in Figure \ref{fig:d}. $\lam(G^{(m)})$ denotes algebraic conflict.} 
			\label{fig:ac}
		\end{figure}

\subsection{Real signed networks and signed corona graphs}

As mentioned in the introduction, we now explore how signed corona graphs can gain us to preserve properties of real signed networks. We consider all the data sets of real signed networks, see Table \ref{tab:data}, that are reported in \cite{leskovec2010signed} and \cite{derr2018signed}, and available in SNAP project webpage \cite{snap}. One remarkable characteristic about these real networks is that each data set contains very high percentage of positive links than negative links, and balance triads dominates in number in all these networks. In particular, the number of $T_0$ is almost $70-90\%$ among all other type of triads in all the networks.

		\begin{table}[H] 
			\begin{center}
			   \begin{tabular}{|c|r|r|r|r|r|}
					\hline
					\textbf{Network} & Nodes & Edges & $p(E_+)$ & Triad &  $p(T)$ \\
					\hline\hline
					\multirow{4}{*}{Epinions} & \multirow{4}{*}{119,217} & \multirow{4}{*}{814,200} & \multirow{4}{*}{0.85} & $T_0$ & 0.870\\
					& & & & $T_1$ &  0.052\\
					& & & & $T_2$ &  0.071\\
					& & & & $T_3$ &  0.007\\
					\hline
					\multirow{4}{*}{Slashdot} & \multirow{4}{*}{82,114} & \multirow{4}{*}{549,202} & \multirow{4}{*}{0.77} &$T_0$ &  0.840\\
					& & & & $T_1$ &  0.077\\
					& & & & $T_2$ &  0.072\\
					& & & & $T_3$ &  0.011\\
					\hline
					\multirow{4}{*}{wiki-Elec} & \multirow{4}{*}{7,118} & \multirow{4}{*}{103,747} & \multirow{4}{*}{0.78} & $T_0$ & 0.702\\
					& & & & $T_1$ &  0.080\\
					& & & & $T_2$ &  0.207\\
					& & & & $T_3$ &  0.011\\
\hline
					\multirow{4}{*}{Bitcoin-OTC} & \multirow{4}{*}{5,881} & \multirow{4}{*}{35,592} & \multirow{4}{*}{0.89} & $T_0$ & 0.724 \\
					& & & & {$T_1$} & 0.122 \\
					& & & & {$T_2$} & 0.145 \\
					& & & & {$T_3$} & 0.009 \\
					\hline 
					\multirow{4}{*}{Bitcoin-Alpha} & \multirow{4}{*}{3,783} & \multirow{4}{*}{24,186} & \multirow{4}{*}{0.93} & $T_0$& 0.793 \\				
					& & & & {$T_1$} & 0.134 \\
					& & & & {$T_2$} & 0.069 \\
					& & & & {$T_3$} & 0.004 \\	
					\hline
				\end{tabular}
			\end{center}
\caption{Data sets of real signed networks. Epinions: A signed social network, Slashdot: Slashdot Zoo signed social network from February 21 2009, wiki-Elec: Wikipedia adminship election data, Bitcoin OTC: Bitcoin OTC web of trust network, Bitcoin Alpha: Bitcoin Alpha web of trust network. The original data sets are signed and directed/weighted, but we consider the undirected/unweighted version of the same data. $p(E_+)$ denotes the proportion of positive edges in the entire network. $p(T)$ is the fraction of triads of type $T$ in the entire network. }
\label{tab:data}
		\end{table}

It can be observed from the the expression of number of each type of triads and signed edges in $G^{(m)}$ obtained in Section \ref{sec:3.1} that the distribution of $T_j$s and signed edges in $G^{(m)}$ are highly influenced by number of $T_j$, signed links, number of nodes and marking of nodes in the seed graph. Besides, the number of nodes in $G^{(m)}$ only can be of the form $n(n+1)^m$ where $n$ denotes the number of nodes in the seed graph. Hence it is evident that $G^{(m)}$ can not have any desired number of nodes in it, but given a total number of nodes in a synthetic network to be produced, one can determine suitable choices of $n, m$ such that the number of nodes in $G^{(m)}$ approximately matches with the desired number of nodes. Then the goal will be to choose a seed graph on $n$ number of nodes that can approximately catch up the desired number of signed edges and triads $T_j$ in $G^{(m)}.$ An another guide from Theorem \ref{Thm:b} is that the marking of nodes which are adjacent to a positive/negative link in the seed graph must follow some conditions in order to produce all types of signed triads in the network.   

Structural properties of corona graphs corresponding to some seed graphs on $4$ or $5$ nodes for which all type of triads exist, are described in Table \ref{tab:sample2}. It follows that these corona graphs fail to preserve properties of real networks. Evidently, a corona graph $G^{(m)}$ for some $m\geq 1$ can preserve structural properties of real networks only when an appropriate seed graph should be chosen. Note that in real networks $\T_0\gg \T_j, j=1,2,3$ and $e_+\gg e_-.$ Thus seed graph must contain high number of positive edges and $T_0$, and negative edges in the seed graph must be connected in such a way that the generation of $T_j, j=1,2,3$ can be controlled during the formation of resulting corona graph.

	\begin{table}[H]
			\begin{center}
				\begin{tabular}{|c|c|r|r|r|r|r|}
					\hline
					$m$ & $G^{(0)}$ & $N$ &  $E$& $p(E_+)$ & Triad  & $p(T_l)$  \\								\hline	
					\hline
					\multirow{4}{*}{6} &\multirow{4}{*}{\begin{tikzpicture}[scale=0.5]
						\draw [fill] (0, 0) circle [radius=0.1];
						\draw [fill] (1, 0) circle [radius=0.1];
						\draw [fill] (0,1) circle [radius=0.1];
						\draw [fill] (1,1) circle [radius=0.1];
						\draw (1,0) --(1,1);
						\draw [dashed] (0,0)	-- (1,0);
						\draw (0,0) -- (0,1);
						\draw (0,1) -- (1,1);
						\end{tikzpicture}} & \multirow{4}{*}{62,500} & \multirow{4}{*}{1,24,996} & \multirow{4}{*}{0.625} & $T_0$  & 0.125  \\
					& & & &  & $T_1$ &  0.625\\	
					& & & &  & $T_2$ &  0.125\\	
					& & & &  & $T_3$ &  0.125\\	
					
					\hline
					\multirow{4}{*}{6} & \multirow{4}{*}{\begin{tikzpicture}[scale=0.5]
						\draw [fill] (0, 0) circle [radius=0.1];
						\draw [fill] (1, 0) circle [radius=0.1];
						\draw [fill] (0,1) circle [radius=0.1];
						\draw [fill] (1,1) circle [radius=0.1];
						\draw [dashed] (1,0) --(1,1);
						\draw (0,0)	-- (1,0);
						\draw (0,0) -- (0,1);
						\draw (0,1) -- (1,1);
						\draw (0,0) -- (1,1);
						\draw (1,0) -- (0,1);
						\end{tikzpicture}} & \multirow{4}{*}{62,500}    & \multirow{4}{*}{1,56,246} & \multirow{4}{*}{0.700} & $T_0$ &  0.250\\
					& & & &  & $T_1$  & 0.650 \\	
					& & & &  & $T_2$  & 0.049\\	
					& & & &  & $T_3$ & 0.049\\			
						\hline
						\hline
					\multirow{4}{*}{5} & \multirow{4}{*}{\begin{tikzpicture}[scale=0.5]
						\draw [fill] (0, 0) circle [radius=0.1];
						\draw [fill] (1, 0) circle [radius=0.1];
						\draw [fill] (0,1) circle [radius=0.1];
						\draw [fill] (1,1) circle [radius=0.1];
						\draw [fill] (.5,1.5) circle [radius=0.1];
						\draw (1,0) --(1,1);
						\draw [dashed] (0,0)	-- (1,0);
						\draw (0,0) -- (0,1);
						\draw (0,1) -- (1,1);
						\draw (.5,1.5) -- (0,0);
						\draw (.5,1.5) -- (0,1) ;
						\draw (.5,1.5) -- (1,1);
						\draw (.5,1.5) --  (1,0);
						\end{tikzpicture}} & \multirow{4}{*}{38,800}   & \multirow{4}{*}{1,01,083} & \multirow{4}{*}{0.754} & $T_0$ &  0.450\\
					& & & & &  $T_1$ &  0.433 \\	
					& & & & &  $T_2$ &  0.050\\	
					& & & & &  $T_3$ &  0.066\\			
					\hline
					\multirow{4}{*}{5} & \multirow{4}{*}{\begin{tikzpicture}[scale=0.5]
						\draw [fill] (0, 0) circle [radius=0.1];
						\draw [fill] (1, 0) circle [radius=0.1];
						\draw [fill] (0,1) circle [radius=0.1];
						\draw [fill] (1,1) circle [radius=0.1];
						\draw [fill] (.5,1.5) circle [radius=0.1];
						\draw (1,0) --(1,1);
						\draw [dashed] (0,0)	-- (1,0);
						\draw (0,0) -- (0,1);
						\draw (0,1) -- (1,1);
						\draw (0,0) -- (1,1);
						\draw (1,0) -- (0,1);
						\draw (.5,1.5) -- (0,0);
						\draw (.5,1.5) -- (0,1) ;
						\draw (.5,1.5) -- (1,1);
						\draw (.5,1.5) --  (1,0);
						\end{tikzpicture}} & \multirow{4}{*}{38,800}   & \multirow{4}{*}{1,16,653} & \multirow{4}{*}{0.786} & $T_0$ &  0.470\\
					& & & &  & $T_1$ &  0.460 \\	
					& & & &  & $T_2$ &  0.030\\	
					& & & &  & $T_3$ &  0.039\\			
					\hline
					\multirow{4}{*}{5} & \multirow{4}{*}{\begin{tikzpicture}[scale=0.5]
						\draw [fill] (0, 0) circle [radius=0.1];
						\draw [fill] (1, 0) circle [radius=0.1];
						\draw [fill] (0,1) circle [radius=0.1];
						\draw [fill] (1,1) circle [radius=0.1];
						\draw [fill] (.5,1.5) circle [radius=0.1];
						\draw [dashed] (1,0) --(1,1);
						\draw [dashed] (0,0) -- (1,0);
						\draw (0,0) -- (0,1);
						\draw (0,1) -- (1,1);
						\draw [dashed] (0,0) -- (1,1);
						\draw (1,0) -- (0,1);
						\draw (.5,1.5) -- (0,0);
						\draw (.5,1.5) -- (0,1) ;
						\draw [dashed] (.5,1.5) -- (1,1);
						\draw (.5,1.5) --  (1,0);
						\end{tikzpicture}} & \multirow{4}{*}{38,800}   & \multirow{4}{*}{1,16,635} & \multirow{4}{*}{0.586} & $T_0$ &  0.180\\
					& & & & &  $T_1$ &  0.500 \\	
					& & & & &  $T_2$ &  0.220\\	
					& & & & &  $T_3$ &  0.100\\			
					\hline
				\end{tabular}
\caption{Statistics of signed links and signed triads of Corona graphs $G^{(m)}$ defined by seed graph $G^{(0)}.$  $N$ and $E$ denote the total number of nodes and links in $G^{(m)}.$ $p(E_+)$ and $p(T_j)$ denote the fraction of positive links and triads of type $T_j$ in $G^{(m)}.$}
\label{tab:sample2}	
			\end{center}				
		\end{table}

Recall that the total number of triads $\T_j$ of type $T_j, j=0,1,2,3$ in $G^{(m)}$ generated by a seed graph $G=G^{(0)}$ is given by \beano  \T_0 &=& \T^{(0)}_0  k_1+ |E^{(0)}_+|^{\stackrel ++} k_2 + \left(|E^{(0)}_+|^{\stackrel --} - |E^{(0)}_+|^{\stackrel ++}\right)k_3\\
 \T_1 &=& \T^{(0)}_1 k_1 + \left( |E^{(0)}_+|^{\stackrel +-} +  |E^{(0)}_-|^{\stackrel ++} \right) k_2 + \left(  |E^{(0)}_-|^{\stackrel --} -  |E^{(0)}_-|^{\stackrel ++} \right) k_3 \\
\T_2 &=&  \T^{(0)}_2k_1 +\left( |E^{(0)}_+|^{\stackrel --} +  \, |E^{(0)}_-|^{\stackrel +-} \right) k_2  +  \left( |E^{(0)}_+|^{\stackrel ++} -  \, |E^{(0)}_+|^{\stackrel --} \right) k_3 \\
\T_3 &=&   \T^{(0)}_3 k_1 + |E^{(0)}_-|^{\stackrel --} k_2 + \left(  |E^{(0)}_-|^{\stackrel ++} -  |E^{(0)}_-|^{\stackrel --} \right) k_3,
\eeano  and the total number of positive and negative edges in $G^{(m)}$ are given by  \beano e_+ &=&  e^{(0)}_+ k_1 + n^{(0)}_+k_2 + (n_-^{(0)} - n_+^{(0)}) k_3 \\ e_- &=& e^{(0)}_- k_1 + n^{(0)}_- k_2 + (n_+^{(0)} - n_-^{(0)}) k_3 \eeano 
respectively, where $k_1=(n+1)^{m}, k_2=(n+1)^m -1, k_3= \sum_{i=0}^{m-1} n_-^{(i)},$ the values of $n_-^{(i)}$ can be found from Table \ref{tab:nmstat}. $\T_j^{(0)}$ denotes the number of triads of type $T_j$ in $G.$ $n_-^{(0)}$ and $n_+^{(0)}$ denote the number of negatively and positively marked nodes  in $G$, and  $e^{(0)}_-$ and $e^{(0)}_+$ denote the number of negative and positive edges in $G$ respectively. However there are six types of edges in the seed graph based on the marking of the adjacent nodes. In order to obtain a compact expression of $\T_j$s in the corona graph we set some of the values of the parameters to be zero that is we exclude certain types of edges in $G.$ This helps to set values of the other parameters so that finally the desired distribution of $\T_j$ can be obtained. 

	\begin{table}[H]
			\begin{center}
				\begin{tabular}{|c|c|c|r|r|r|r|r|}
					\hline
					Graph & $m$ & $G^{(0)}$ & $N$ & $E$ &  $p(E_+)$ & Triad  & $p(T_l)$  \\														
					\hline \hline
					\multirow{4}{*}{1} &\multirow{4}{*}{3} & \multirow{4}{*}{\begin{tikzpicture}[scale=0.5]
						\draw [fill] (0, 0) circle [radius=0.1];
				\draw [fill] (1, 0) circle [radius=0.1];
				\draw [fill] (3,0) circle [radius=0.1];
				\draw [fill] (0,1) circle [radius=0.1];
                                          \draw [fill] (1,1) circle [radius=0.1];
                                           \draw [fill] (2,1) circle [radius=0.1];
                                           \draw [fill] (0,-1) circle [radius=0.1];
                                           \draw [fill] (1,-1) circle [radius=0.1];
                                           \draw [fill] (3,-1) circle [radius=0.1];
				\draw [dashed] (0,0) --(0,1);
				\draw (0,0) --(1,0);
				\draw (0,0) --(1,-1);
				\draw (1,0) --(3,0);
				\draw (0,0) --(2,1);
				\draw (1,0) --(2,1);
				\draw (1,0) --(1,-1);
				\draw (1,0) --(3,-1);
				\draw (1,-1) --(3,-1);
				\draw (1,-1) --(2,1);
				\draw (1,-1) --(3,0);
				\draw (3,-1) --(2,1);
				\draw (3,-1) --(3,0);
				\draw (2,1) --(3,0);
				\draw [dashed] (1,0)-- (0,-1);
				\draw [dashed] (1,0)-- (1,-1);
				\draw [dashed] (0,0)-- (1,1);
				\draw [dashed] (1,0)-- (1,1);
						\end{tikzpicture}} & \multirow{4}{*}{9,000}   & \multirow{4}{*}{25,991}  & \multirow{4}{*}{0.734} & $T_0$  & 0.762\\
					&& & & & & $T_1$ &  0.053 \\	
					&& & & & & $T_2$ &  0.173\\	
					& && & & & $T_3$ &  0.012\\			
					\hline
\multirow{4}{*}{2} &\multirow{4}{*}{3} &\multirow{4}{*}{\begin{tikzpicture}[scale=0.5]
						\draw [fill] (0, 0) circle [radius=0.1];
				\draw [fill] (1, 0) circle [radius=0.1];
				\draw [fill] (-1,0) circle [radius=0.1];
				\draw [fill] (2,0) circle [radius=0.1];
                                          \draw [fill] (-1,-1) circle [radius=0.1];
                                           \draw [fill] (-1,-2) circle [radius=0.1];
                                           \draw [fill] (0,-3) circle [radius=0.1];
                                           \draw [fill] (1,-2) circle [radius=0.1];
				\draw [dashed] (0,0) --(1,0);
				\draw [dashed] (0,0) --(-1,0);
				\draw [dashed] (1,0) --(2,0);
				\draw (0,0) --(-1,-1);
				\draw (0,0) --(1,-2);
                                          \draw (-1,-1) --(-1,-2);
                                          \draw (-1,-1) --(0,-3);
				\draw (-1,-1) --(1,-2);
				\draw (-1,-2) --(1,-2);
				\draw (0,0) --(0,-3);
				\draw (0,0) --(-1,-2);
                                          \draw (0,-3) --(1,-2);
				\draw (-1,-2) --(0,-3);
				\draw (1,0) --(1,-2);
				\draw (1,0) --(0,-3);
						\end{tikzpicture}} & \multirow{4}{*}{5,832} & \multirow{4}{*}{16,759} & \multirow{4}{*}{0.739} & $T_0$ &  0.714  \\
					& && & & & $T_1$ &  0.098\\	
					& && & & & $T_2$ &  0.178\\	
					& & && & & $T_3$  & 0.009\\
\hline
\multirow{4}{*}{3} &\multirow{4}{*}{4} &\multirow{4}{*}{\begin{tikzpicture}[scale=0.5]
						\draw [fill] (0, 0) circle [radius=0.1];
				\draw [fill] (1, 0) circle [radius=0.1];
				\draw [fill] (-1,0) circle [radius=0.1];
				\draw [fill] (2,0) circle [radius=0.1];
                                          \draw [fill] (-1,-1) circle [radius=0.1];
                                           \draw [fill] (-1,-2) circle [radius=0.1];
                                           \draw [fill] (0,-3) circle [radius=0.1];
                                           \draw [fill] (1,-2) circle [radius=0.1];
                                          \draw [fill] (2,-2) circle [radius=0.1];
				\draw [dashed] (0,0) --(1,0);
				\draw [dashed] (0,0) --(-1,0);
				\draw [dashed] (1,0) --(2,0);
				\draw (0,0) --(-1,-1);
				\draw (0,0) --(1,-2);
                                          \draw (-1,-1) --(-1,-2);
                                          \draw (-1,-1) --(0,-3);
				\draw (-1,-1) --(1,-2);
				\draw (-1,-2) --(1,-2);
				\draw (0,0) --(0,-3);
				\draw (0,0) --(-1,-2);
                                          \draw (0,-3) --(1,-2);
				\draw (-1,-2) --(0,-3);
				\draw (1,0) --(1,-2);
				\draw (1,0) --(0,-3);
				\draw [dashed] (2,-2) --(2,0);
				\draw (1,-2) --(2,0);
						\end{tikzpicture}} & \multirow{4}{*}{90,000} & \multirow{4}{*}{2,59,991} & \multirow{4}{*}{0.738} & $T_0$ &  0.705 \\
					&& & & & & $T_1$ &  0.151\\	
					&& & & & & $T_2$ &  0.133\\	
					& && & & & $T_3$  & 0.010\\
\hline
				\end{tabular}
\caption{Corona graph $G^{(m)}$ which preserves characteristics of signed links and triads of real networks. $N$ and $E$ denote the total number of nodes and links in $G^{(m)}.$ $p(E_+)$ and $p(T_j)$ denote the fraction of positive links and triads of type $T_j$ in $G^{(m)}$ respectively.}
\label{tab:new2}
			\end{center}				
		\end{table}

Setting $|E^{(0)}_+|^{\stackrel +-}=|E^{(0)}_+|^{\stackrel --}=|E^{(0)}_-|^{\stackrel --}=0$ and $\T^{(0)}_3 =0,$ we obtain  \beano  \T_0 &=& \T^{(0)}_0  k_1+ |E^{(0)}_+|^{\stackrel ++} (k_2-k_3)\\
 \T_1 &=& \T^{(0)}_1 k_1 + |E^{(0)}_-|^{\stackrel ++} (k_2-k_3) \\
\T_2 &=&  \T^{(0)}_2k_1 + |E^{(0)}_-|^{\stackrel +-} k_2 +  |E^{(0)}_+|^{\stackrel ++}k_3 \\
\T_3 &=&    |E^{(0)}_-|^{\stackrel ++} k_3.\eeano
The reason behind choosing $\T^{(0)}_3=0$ is that we want the number of $T_3$ to be minimum among all triads. Besides, note that $k_2\gg k_3.$ Finally considering $e_+^{(0)} \gg e_-^{(0)}, \T^{(0)}_0\gg \T_j^{(0)}, j=1,2$ the desired statistics of real networks can be obtaind. In Table \ref{tab:new2} we consider certain seed graphs under these conditions and show that the corresponding corona graphs can preserve the signed edge and triad distribution of real networks. The positive and negative degree distributions of the corona graphs $G^{(m)}$ defined by the seed graphs given in Table \ref{tab:new2} are plotted in Figure \ref{fig:pn_deg_dist2}. However do not compare it with real networks since the number of nodes and links in $G^{(m)}$ do not match the same for real networks.  

\begin{figure}[H]
			\begin{subfigure}[b]{0.5\textwidth}
				\centering
				\includegraphics[scale=0.4]{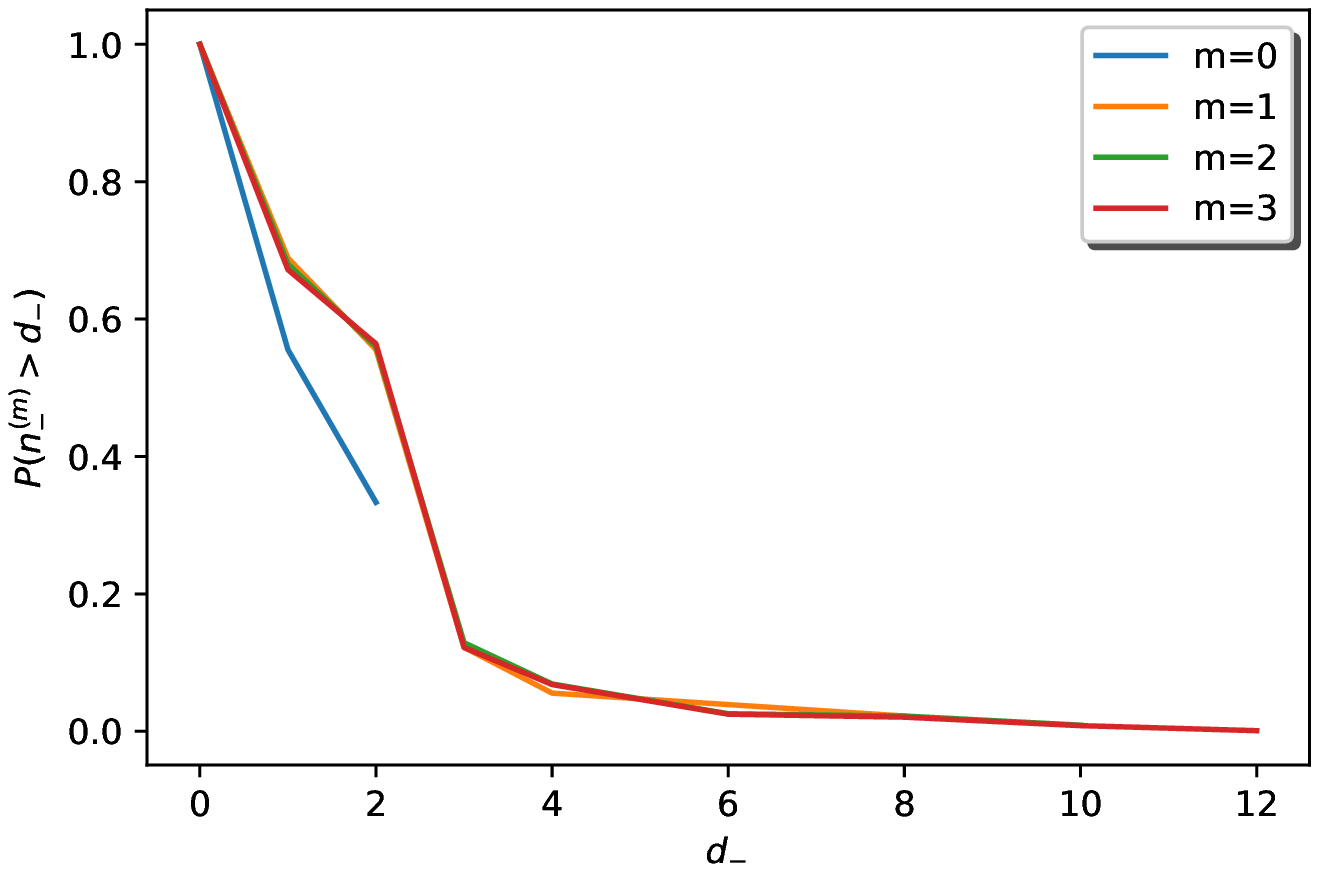}
				\caption{Negative degree distribution of $G^{m}$ generated by Graph $1$}
			\end{subfigure}
			\begin{subfigure}[b]{0.5\linewidth}
				\centering
				\includegraphics[scale=0.4]{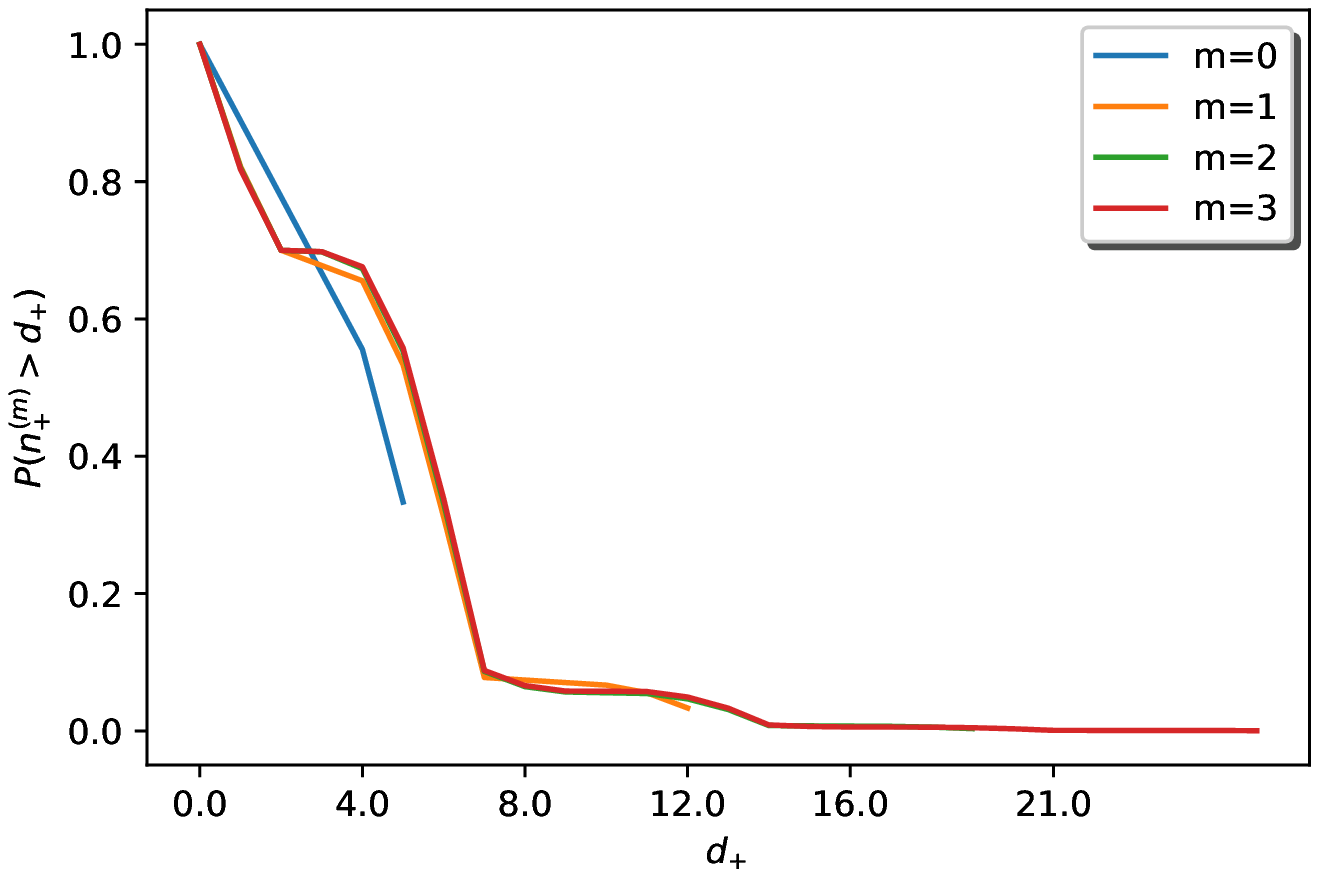}
				\caption{Positive degree distribution of $G^{m}$ generated by Graph $1$}
			\end{subfigure}
			\begin{subfigure}[b]{0.5\linewidth}
				\centering
				\includegraphics[scale=0.4]{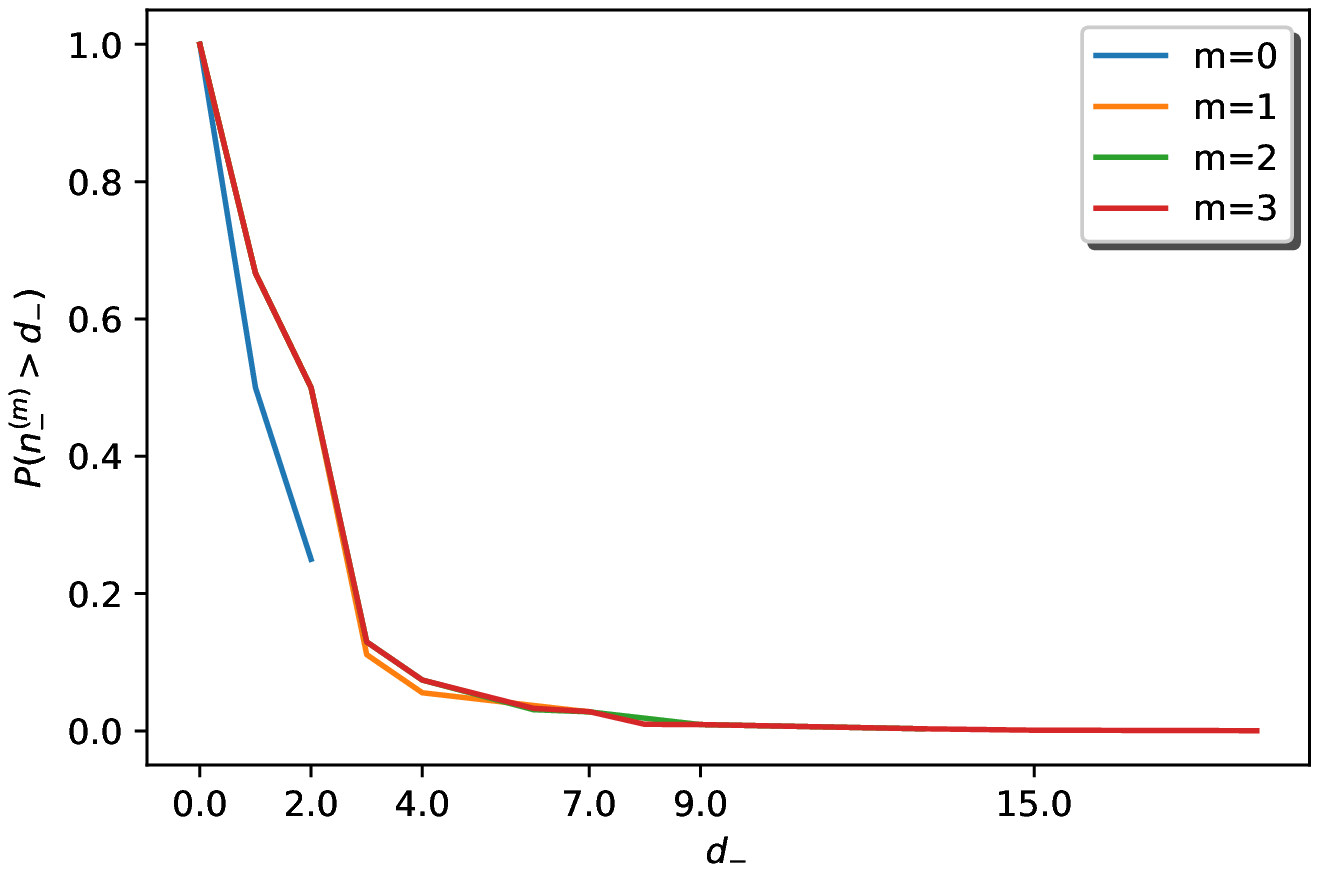}
				\caption{Negative degree distribution of $G^{m}$ generated by Graph $2$}
			\end{subfigure}
			\begin{subfigure}[b]{0.5\linewidth}
				\centering
				\includegraphics[scale=0.4]{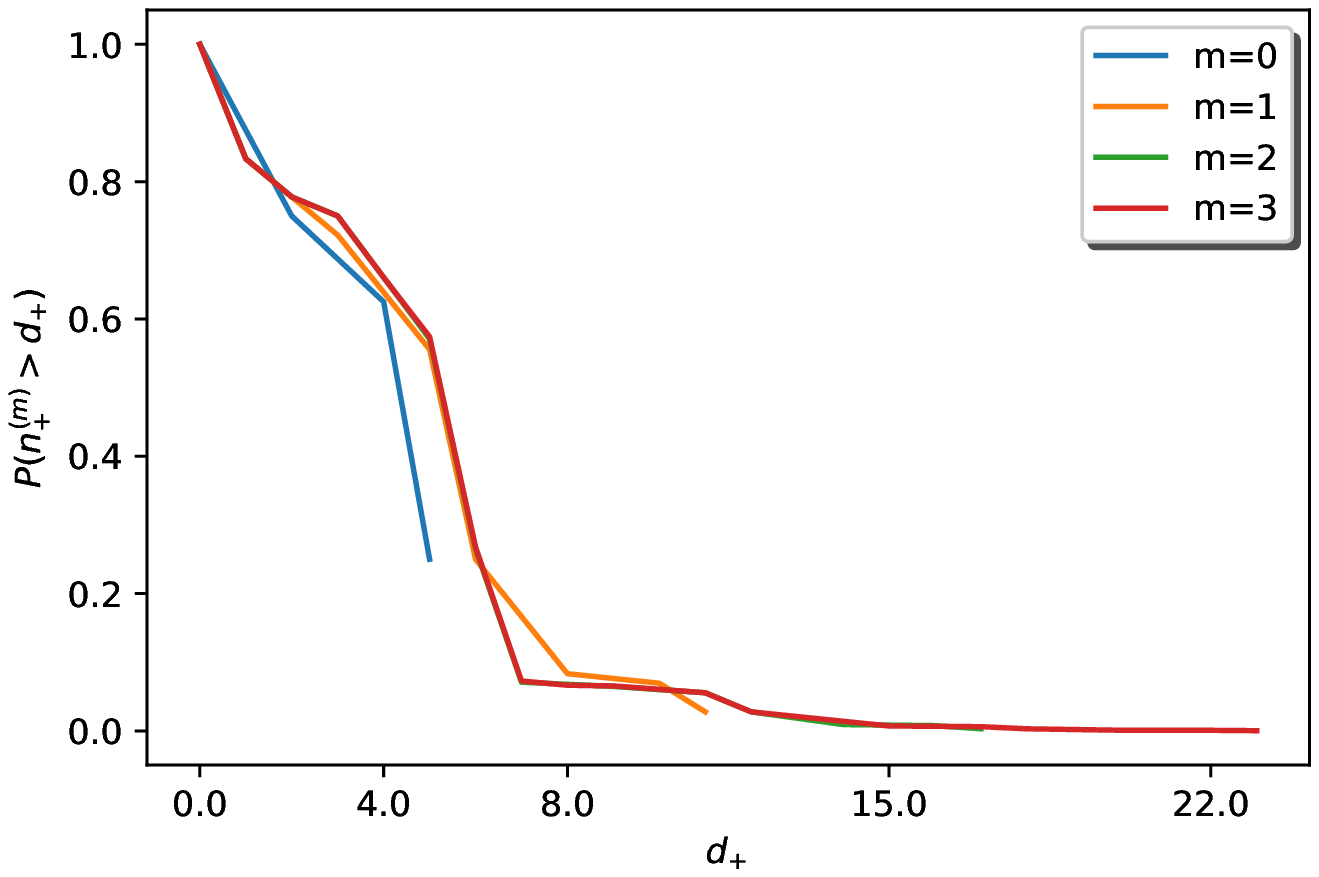}
				\caption{Positive degree distribution of $G^{m}$ generated by Graph $2$}
                                \end{subfigure}
                                 \begin{subfigure}[b]{0.5\linewidth}
				\centering
				\includegraphics[scale=0.4]{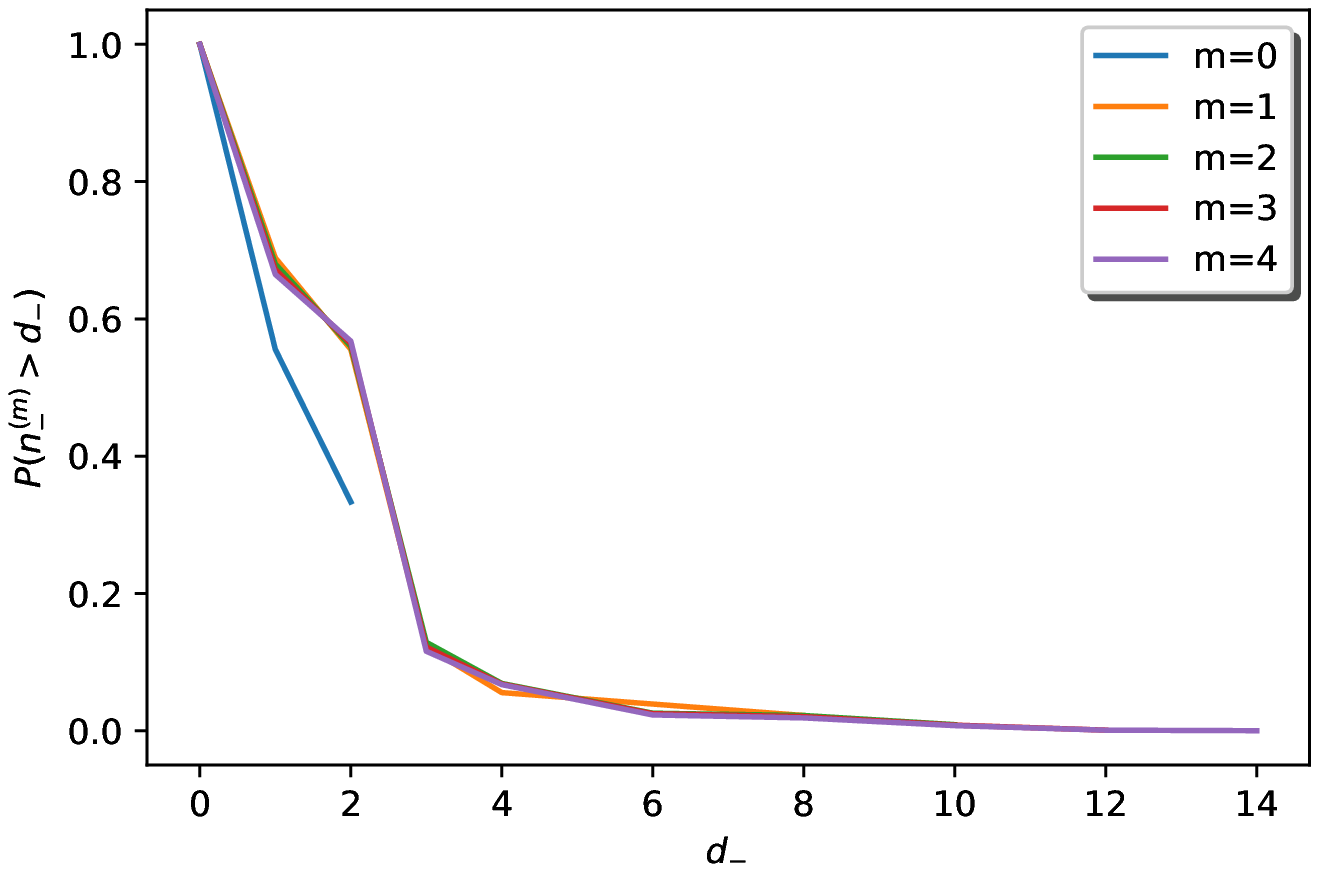}
				\caption{Negative degree distribution of $G^{m}$ generated by Graph $3$}
			\end{subfigure}
			\begin{subfigure}[b]{0.5\linewidth}
				\centering
				\includegraphics[scale=0.4]{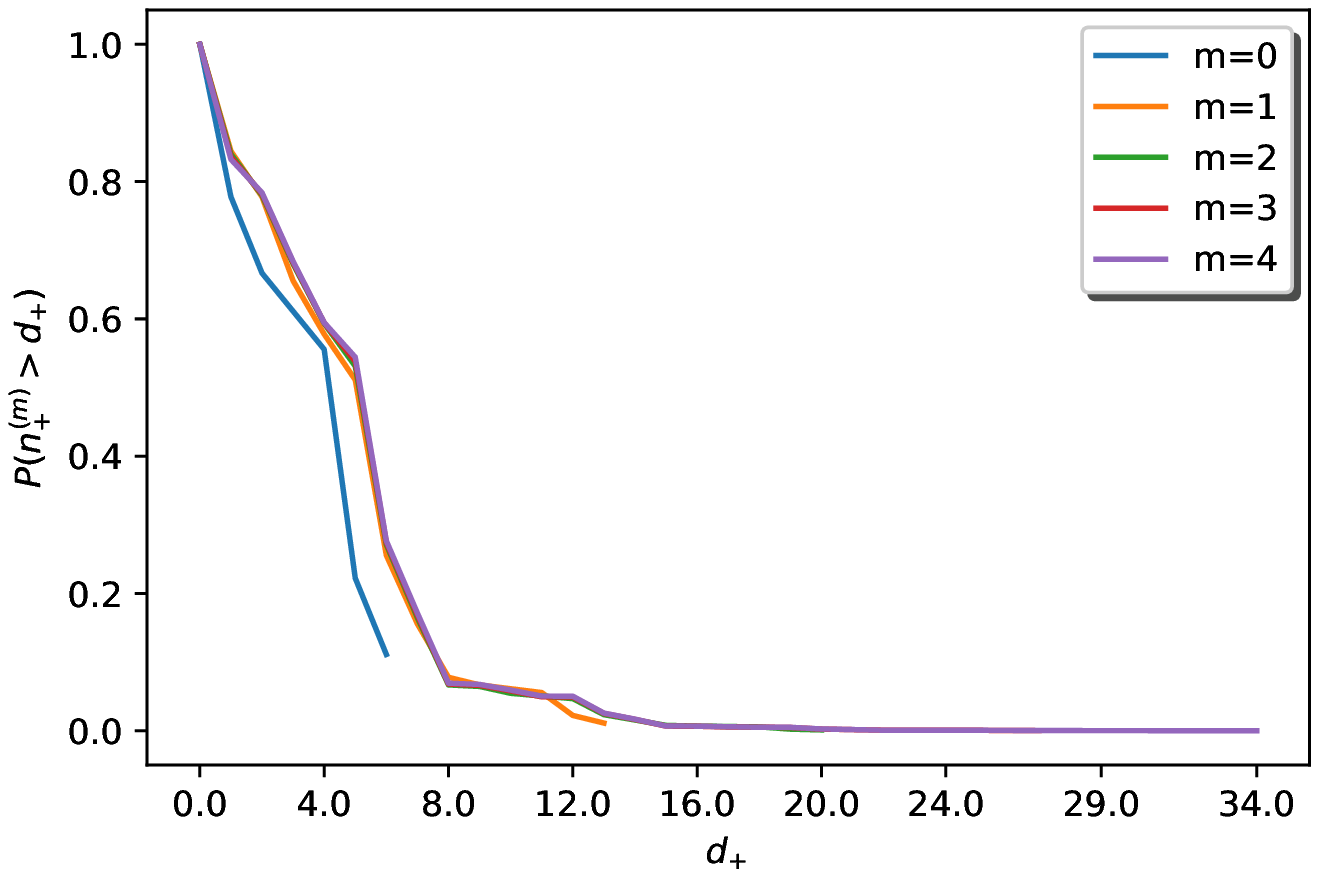}
				\caption{Positive degree distribution of $G^{m}$ generated by Graph $3$}
			\end{subfigure}
			\caption{(a) and (c) ((b) and (d)) shows the fraction of nodes in $G_1^{(m)}$ and  $G_2^{(m)}$ having degree greater than or equal to $d_-$ ($d_+$), respectively.} 
			\label{fig:pn_deg_dist2}
		\end{figure}

\section{Conclusion}

We defined corona product of two signed graphs and study their structural and spectral properties. A signed network generative model is proposed based on corona product by taking the corona product of a small graph with itself iteratively. This small graph is called seed graph for the resulting signed graphs, which are called corona graphs. We derived fundamental properties of corona graphs that include number of signed links, signed triads, degree distribution and algebraic conflict of the corona graphs. Finally we show that a seed graph can be chosen for which the corresponding corona graphs can preserve properties of real signed networks.

We believe that the framework of corona product introduced in this paper can be extended to define corona product of gain graphs and weighted graphs. Besides, generalized corona product of signed graphs can be defined alike generalized corona product for unsigned graphs \cite{laali2016spectra}.



\end{document}